\documentclass[12pt,a4]{elsarticle}

\usepackage{latexsym,amsbsy,amssymb,amsmath}
\usepackage{color,graphicx}
\usepackage{color}
\usepackage{geometry}
\usepackage{epstopdf}

\usepackage[justification=centering]{caption}
\usepackage{subcaption}

\renewcommand{\div}{\operatorname{div}}
\newcommand{\card}{\operatorname{card}}

\newcommand{\conv}{\operatorname{conv}}

\begin{document}
\begin{frontmatter}

\title{On stability, convergence and accuracy of bES-FEM and bFS-FEM for nearly incompressible elasticity}

\author[hai]{Thanh Hai Ong}
\author[clairekye]{Claire E. Heaney}
\author[clairekye]{Chang-Kye Lee}
\author[liu]{G.R. Liu}
\author[hung1,hung2]{H. Nguyen-Xuan\corref{nxh}}
\cortext[nxh]{Corresponding author. Email address: hung.nx@vgu.edu.vn (H. Nguyen-Xuan).}
\address[hai]{Department of Analysis, Faculty of Mathematics Computer Science,University of Science,VNU-HCMC, Nguyen Van Cu Street,  District 5, Ho Chi Minh City, Vietnam}
\address[clairekye]{Institute of Mechanics and Advanced Materials, School of Engineering, Cardiff University, Queen's Buildings, The Parade, Cardiff CF24 3AA, UK}
\address[liu]{School of Aerospace Systems, University of Cincinnati, 2851 Woodside Dr, Cincinnati, OH 45221, USA}
\address[hung1]{Department of Computational Engineering, Vietnamese-German University, Binh Duong New City, Vietnam}
\address[hung2]{Department of Architectural Engineering, Sejong University, 98 Kunja Dong, Kwangjin Ku, Seoul 143-747, South Korea}

\begin{abstract}
We present in this paper a rigorous theoretical framework to show stability, convergence and accuracy of improved edge-based and face-based smoothed finite element methods (bES-FEM and bFS-FEM) for nearly-incompressible elasticity problems. The crucial idea is that the space of piecewise linear polynomials used for the displacements is enriched with bubble functions on each element, while the pressure is a piecewise constant function. The meshes of triangular or tetrahedral elements required by these methods can be generated automatically. The enrichment induces a softening in the bilinear form allowing the weakened weak ($W^2$) procedure to produce a high-quality solution, free from locking and that does not oscillate. We prove theoretically that both methods confirm the uniform inf-sup and convergence conditions. Four numerical examples are given to validate the reliability of the bES-FEM and bFS-FEM.
\end{abstract}
\end{frontmatter}
\section{Introduction} \label{sec;introduction}
Rubber-like materials are able to withstand extremely high strains whilst exhibiting very little or no permanent deformation and consequently are widely used in industry. In addition to elastic properties, the volume of these materials is almost preserved upon loading. Rubber-like materials are said therefore to be nearly incompressible and typically possess bulk moduli that are several orders of magnitude higher than their shear moduli (equivalently, they have a Poisson's ratio close to one half). It is well known that the stress analysis of nearly-incompressible materials requires special care. Applying low-order finite elements based on quadrilaterals, hexahedra, triangles or tetrahedra, to such problems, results in a severe underprediction of the displacement known as locking. A variety of numerical methods have been proposed to overcome this defect, for example: $h$-version of finite elements \cite{BM1, BM2}, B-bar method \cite{H}, mixed
formulations \cite{AW, B1}, enhanced assumed strain (EAS) modes \cite{CYMCAF, SR}, reduced integration stabilization \cite{ FO} and two-field
mixed stress elements \cite{PW}, a stream function approach \cite{FLA} and mimetic finite difference method \cite{LK} and so on. In addition to these, several publications investigate an average nodal pressure formulation in which a constant pressure field is enforced over a patch of triangles or tetrahedra \cite{BoBu, BMH, NPO, LAMI, PNC}. Despite the many available approaches for solving nearly-incompressible elasticity problems on a triangulation, only a few methods are based on rigorous mathematical analysis. An example of one such method can be found in \cite{LAMI}. Here, the author introduced a discontinuous pressure and used bubble functions in order to enrich the space of piecewise linear polynomials to which the displacements belong. However, the method still has certain drawbacks inherited from FEM such as 1) an overestimation of the stiffness matrix for nearly-incompressible and bending-dominated problems, 2) a poor performance for distorted meshes, 3) a poor accuracy of the stresses. Moreover, we make mention of the very important three-field (Hu-Washizu) methods. In fact many of the two-field methods mentioned in the overview can be derived as special cases of the Hu-Wahsizu formulation, for which a rigorous analysis has been carried out in \cite{DLRW,LRW}.\\

In this paper we propose two improved methods which use bubble functions as enrichments to the edge-based and face-based smoothed finite element methods (bES-FEM and bFS-FEM). These methods contribute to the further development of advanced numerical tools that can be used for nearly-incompressible elasticity problems, whilst simultaneously building on the advantages of some classical methods as explained below.\\

Firstly, an improved version of a so-called bES-FEM has the same desirable features as bES-FEM-T3 studied in \cite{HG}. Both bES-FEM and bFS-FEM work well for three-dimension problems, where bubble functions are generally defined by the $(d+1)$th-power bubble function and the hat function. Most importantly, both methods are theoretically proven to ensure the uniform inf-sup condition and the convergence. In addition, there is a basic difference between bES-FEM and bES-FEM-T3 as follows: for bES-FEM, the approximate pressure and displacement are directly computed by the mixed approach provided in (\ref{phtrinh4}a) and (\ref{phtrinh4_pressure}) while for bES-FEM-T3, the approximate pressure is computed as a posteriori of the displacements based on the edge-based smoothing domains.\\

Secondly, we use mixed methods \cite{AW, BF} to reformulate the linear elasticity problem as a mixed displacement-pressure problem. Our aim is to attain a good approximation to the pressure solution \cite{BM1}, which we model here as piecewise constant.\\

Thirdly, the proposed approximation to the displacement solution is a combination of the displacement from ES-FEM/FS-FEM \cite{LTK, TLKG} and the displacement from the bubble functions \cite{RP, JM}. ES-FEM and FS-FEM improved the standard FE strain fields via a strain smoothing technique described in \cite{CWYY}.  The methods proposed in this paper build on ES-FEM and FS-FEM, and therefore inherit the positive qualities associated with this smoothing technique, namely 1) its solutions are more accurate than those of linear triangular elements (FEM-T3) and quadrilateral elements (FEM-Q4) using the same sets of nodes \cite{H, ZT}; 2) ES-FEM and FS-FEM perform well with distorted meshes; 3) their stress solutions, which are very precise, have the convergent property and 4) they can be easily implemented into existing FEM packages without requiring additional degrees of freedom. Clearly this technique of smoothing is a powerful tool and it has already been applied to a wide range of practical mechanics problems, e.g, \cite{LNZ, HRBD, HT}. Nevertheless, if the displacement is approximated only by ES-FEM or FS-FEM {\it i.e.} without enrichment by bubble functions, these methods violate the inf-sup condition and uniform convergence. Other methods from the SFEM family also fail to satisfy this condition, implying that they also suffer from volumetric locking in the case of nearly-incompressible elasticity \cite{TLKG, HG}. To overcome volumetric locking for the SFEM family, only a few approaches have been presented. For example, in \cite{TLKG}, the authors suggested a combined FS/NS-FEM model and in \cite{HG} the use of bubble functions was proposed. Neither of these approaches are based on a rigorous mathematical analysis.\\

Finally, the degree of freedom which is associated with the pressure variable can be statically condensed out of the system of equations, in contrast to the method based on the classical MINI element \cite{BF}, for example, where condensation cannot be applied.\\

The rest of this paper is organized as follows. In the next section, we briefly recall the boundary value problem of linear elasticity, the mixed displacement-pressure formulation and its associated weak form. Section~\ref{sec:bES-bFS} describes the enrichment of ES-FEM and FS-FEM by bubble functions. Section~\ref{mathproper} presents the mathematical properties of bES-FEM and bFS-FEM, where only small deformations are considered. Displacement, energy and pressure error norms are defined in section~\ref{Error} for the precise quantitative examination of various models. Four numerical tests are presented in section~\ref{Num} to demonstrate the effectiveness and accuracy of the proposed methods. In the final test we apply the proposed bES-FEM to a large deformation problem. In the last section we draw conclusions and give possible directions for future work.

\section{The boundary value problem of linear elasticity} \label{sec:bvp_le}
We consider a static linear elasticity problem in a bounded domain $\Omega \subset \mathbb R^d$ , $d = \{2, 3\}$ with a Lipschitz boundary $\partial \Omega$. The governing equations express equilibrium between the Cauchy stresses $\boldsymbol{\sigma}$ and the applied body forces $\bold{f}$
\begin{equation}
\label{phtrinh01}
-\div \boldsymbol{\sigma}  = \bold{f} \quad \text{in}~\Omega.
\end{equation}
The displacement $\bold{u}$ is prescribed on the boundary $\partial \Omega$ by
\begin{equation}
\label{phtrinh02}
\bold{u} = 0 \quad \text{on}~\partial \Omega.
\end{equation}
In addition to (\ref{phtrinh01}) and (\ref{phtrinh02}), we introduce the infinitesimal strain tensor $\boldsymbol{\varepsilon}$ which is related to the displacement $\bold{u}$ by
\begin{equation}
\label{phtrinh03}
\boldsymbol{\varepsilon}_{ij}(\bold{u}) = \frac{1}{2}\left({\partial_j \bold{u}_i + \partial_i \bold{u}_j} \right)  ~~ \text{in}~\Omega, ~~ \forall i,j = \overline{1,d},
\end{equation}
where $\partial_i = \frac{\partial}{\partial x_i}$, $(x_1,\cdots,x_d) \in \mathbb R^d$ and $\boldsymbol{\varepsilon}(\bold{u})=[\boldsymbol{\varepsilon}_{ij}(\bold{u})]_{{i,j}=\overline{1,d}}$. For an isotropic linear elastic material, the constitutive relation is given by
\begin{equation}
\label{phtrinh04}
\boldsymbol \sigma_{ij}(\bold{u}) = \lambda \delta_{ij} \boldsymbol{\varepsilon}_{kk}(\bold{u}) + 2 \mu \boldsymbol{\varepsilon}_{ij}(\bold{u}) \quad \text{in}~\Omega 
\end{equation}
where $\lambda$ and $\mu$ are the Lam\'e constants and $\delta_{ij}$ is the Kronecker delta. 
The Lam\'e constants are related to the Young's modulus, $E$, and Poisson's ratio, $\nu$, through the following:
\begin{equation*}
\lambda = \frac{\nu E}{(1+\nu)(1-2\nu)}, \qquad \mu = \frac{E}{2(1+\nu)}\ .
\end{equation*}

In this paper our attention is devoted to the study of nearly-incompressible materials for which Poisson's ratio is close to $0.5$. Such a choice of this parameter is well known to lead to a poor performance by FEM due to locking and instability.
\subsection{Mixed displacement-pressure formulation and the weak form}
The elasticity problem (\ref{phtrinh01}) can be rewritten in a mixed displacement-pressure form
\begin{eqnarray}
\label{phtrinh04a}
-\div~\boldsymbol{\sigma} &=& \bold{f} \quad \text{in}~\Omega,\\
 \div~\bold{u} - \frac{p}{\lambda} &=& 0 \quad \text{in}~\Omega, \label{phtrinh04_pressure}
\end{eqnarray}
where the pressure $p$ is introduced as an additional variable. The mixed form is equivalent to the penalized Stokes equations. We now introduce several function spaces which are required for the weak form:

\begin{displaymath}
V_0  =  [ H_0^1(\Omega)]^d,    \quad
	L^2_0(\Omega):  =  \bigg\{ q \in L^2(\Omega) : \int\limits_\Omega  q~\text{d}\Omega  = 0  \bigg\}\ .
\end{displaymath}
The space to which the pressure solution belongs is $L^2_0(\Omega)$. The condition that the volume integral of the pressure should be zero follows directly from  integrating equation~\eqref{phtrinh04_pressure}, transforming the integral to a boundary integral and then using the fact that the displacement satisfies homogeneous Dirichlet boundary conditions.
 The mixed approach aims to find a displacement field $\bold{u} \in V_0 $ and a pressure $p \in L^2_0(\Omega)$ that satisfy
\begin{subequations}
\begin{align}
\label{phtrinh05}
a(\bold{u,v}) + b(\bold{v},p) &= (\bold{f,v}) &   & \hspace{-2cm} \forall \bold{v} \in V_0,\\
\label{phtrinh05.1}
b(\bold{u},q) - \frac{1}{\lambda}c(p,q) &= 0  &   & \hspace{-2cm} \forall q\in L^2_0(\Omega).
\end{align}
\end{subequations}
{The bilinear forms are defined as follows:
\begin{eqnarray*}
a(\bold{u,v}) &=& 2\mu \int\limits_\Omega  {\boldsymbol{\varepsilon} ^T (\bold{v}(\bold{x}))\,\bold{D}\,\boldsymbol{\varepsilon} (  \bold{u}(\bold{x}))~\text{d}\Omega}, \quad
b(\bold{u},q) = \int\limits_\Omega  {q(\bold{x})\nabla \cdot \bold{u}(\bold{x})~\text{d}\Omega}, \\
c(q,p) &=& \int\limits_\Omega  {q(\bold{x})p(\bold{x})~\text{d}\Omega}, \quad
(\bold{f,v}) = \int\limits_\Omega { \bold{v}^T (\bold{x}) \, \bold{f}(\bold{x})~\text{d}\Omega},
\end{eqnarray*}}
where $\bold{f} \in [L^2(\Omega)]^d$. In the definition of the bilinear forms we have introduced Voigt notation, in which the components of the stress and strain tensors are arranged in column vectors, for example: \mbox{$\boldsymbol{\varepsilon} = \{\varepsilon_{xx}\ \varepsilon_{yy}\ \varepsilon_{zz}\ \varepsilon_{xy}\ \varepsilon_{yz}\ \varepsilon_{zx}\}^T$}. The matrix $\bold{D}$ of material constants is symmetric, positive definite and its eigenvalues are bounded in $[\lambda^{\bold D}_{\min},\lambda^{\bold D}_{\max}] \subset \mathbb R^+$.
\section{Description of ES-FEM and FS-FEM enriched by bubble functions}
\label{sec:bES-bFS}
\subsection{The finite spaces}

The polygonal domain $\Omega$ is discretized by the triangulation $\mathcal T_h$ (the primal mesh), where $\mathcal T_h$ consists of triangles (2D) or tetrahedra (3D). The set $\mathcal T_h$ has $N_e$ elements, $N_n$ nodes (or vertices), $N_s$ edges, $N_f$ faces (3D) and $\overline \Omega = \bigcup\limits_{i = 1}^{N_e } {\overline{T_i} }$. For each element  $T \in \mathcal T_h$, the barycentric point $c_T$ is called a mesh point of $T$. Let $ \mathcal V_h$ be the standard linear finite element space defined on the triangulation $\mathcal T_h$,
\begin{displaymath}
\mathcal V_h = \left\{ {\bold{u} \in V_0 ,\,\,\bold{u}|_T  \in \left[ {\mathbb P^1 (\Omega )} \right]^d,~\text{for all}~T \in \mathcal T_h } \right\}
\end{displaymath}
which  has the standard nodal basis functions $N_i$ ($i = \overline{1,N_n}$) associated with node $i$. We define the space of bubble functions as
\begin{displaymath}
\mathcal B_h  = \left\{ {b_T  \in \mathbb H^1 (\Omega ),\,\,\left. {b_T } \right|_{\partial T}  = 0\,\,\text{and}\,\,\int\limits_T {b_T(\bold{x}) \text{d}\Omega > 0,T \in \mathcal T_h } } \right\},
\end{displaymath}
where the basis bubble functions are chosen to be one of two types (see \cite{RP} and \cite{JM}).\\

For the first type, the $\xi$th-power bubble function is used for each element $T \in \mathcal T_h$ with $\xi = d+1$
\begin{equation}
\label{cubic}
b_T (\bold{x}) = \left\{ \begin{array}{ll}
 c_b (d+1)^3\prod\limits_{i = 1}^{d+1} {\lambda _{T^{(i)} } (\bold{x})} \quad & \text{if}\,\bold{x} \in T \subset \mathcal T_h,\\
 0                                                                       & \text{elsewhere} \\
 \end{array} \right.
\end{equation}
where each function $\lambda_{T^{(i)}}$ is a barycentric coordinate associated with a vertex $\bold{x}_{T^{(i)}}$ of the triangle $T$, and $c_b$ is computed in such a way that $b_T(c_T) = 1$ where $c_T$ is the centroid of $T$.\\

The second type is a hat function defined on $T$, where $T$ is partitioned into sub-triangles (2D) or sub-tetrahedra (3D), $\{T_{(i)}\}_{i=\overline{1,d+1}}$. This is achieved by joining the centroid $c_T$ to the two vertices on each edge of the triangle in turn (2D), or to the three vertices on each face (3D).
\begin{equation}
\label{hat}
b_T (\bold{x}) = \left\{ \begin{array}{ll}
 c_b (d+1)\lambda _{T_{(i)} } (\bold{x}) \quad & \text{if }\,\bold{x} \in T_{(i)} \subset T\subset\mathbb R^d, \\
 0                                       & \text{elsewhere}. \\
 \end{array} \right.
\end{equation}
The finite element space for the displacement which is enriched with bubble functions is defined as
\begin{displaymath}
\mathcal V^{\mathcal B}_h = \mathcal V_h  \oplus  [\mathcal B_h]^d \subset  [H^1(\Omega)]^d.
\end{displaymath}
Each function $\bold{u}_h \in \mathcal V^{\mathcal B}_h$ which is restricted on $T \in \mathcal T_h$ is written as
\begin{equation}
\label{phtrinh1}
{{\bf{u}}_h}(x) = \underbrace {\sum\limits_{i = 1}^{d + 1} {\left( {{N_{{T^{(i)}}}}(x){\bold{I}}{{\bold{d}}_d}} \right){{\bold{u}}_{{T^{(i)}}}}} }_{{\ell _h}(x) \in {\mathcal V_h}} + \underbrace {\left( {N_{{c_T}}^b(x){\bf{I}}{{\bf{d}}_d}} \right){{\bf{u}}_{{c_T}}}}_{{\bold b_h}(x) \in {{[{\mathcal B_h}]}^d}},\end{equation}
where the identity matrix of size $d$ is denoted by $\bold {Id}_d$, $N_{T^{(i)}}$ is the standard nodal basis function associated with the vertex $\bold{x}_{T^{(i)}}$ of the triangle $T$, $N_{c_T }^b$ is the standard nodal basis bubble function defined on $T$ with the centroid $c_T$. The values $\bold{u}_{T^{(i)}} $ and $\bold{u}_{c_T} \in \mathbb R^d$ are the nodal values of $\bold{u}_h$ at the vertex  $\bold{x}_{T^{(i)}}$ and the barycenter $c_T$.\\
\subsection{The dual mesh}
\begin{figure}[htbp!]
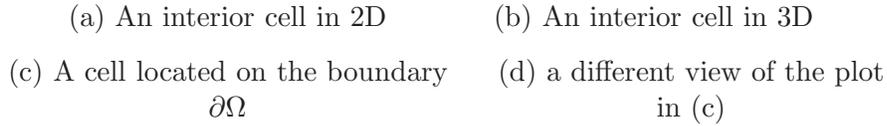

        \centering
        \begin{subfigure}[b]{5cm}
                \caption{An interior cell in 2D}
                \label{fig:cell_bESFEM_a}
        \end{subfigure}%
        \quad 
        \begin{subfigure}[b]{5.5cm}
                \caption{An interior cell in 3D}
                \label{fig:cell_bESFEM_b}
        \end{subfigure}

        \begin{subfigure}[b]{6cm}
                \caption{A cell located on the boundary $\partial \Omega$}
                \label{fig:cell_bESFEM_c}
        \end{subfigure}%
        \quad 
        \begin{subfigure}[b]{5.5cm}
                \caption{a different view of the plot in (c)}
                \label{fig:cell_bESFEM_d}
        \end{subfigure}
        \caption{Illustrations of smoothing cells for bES-FEM.}
        \label{fig:smoothingCell_bESFEM}
\end{figure}

Now, we design a dual mesh for smoothing the strain and the divergence operator. For each of the 2D bES-FEM, the 3D bES-FEM and the 3D bFS-FEM, a dual mesh $\mathcal  T^*_h$ is created in a similar manner to the 2D and the 3D edge-based smoothing domain \cite{LT, CM} and the face-based smoothing domain \cite{LT} respectively. It is constructed  by connecting all vertices, center points of elements in $\mathcal T_h$ and center points of faces (for 3D bES-FEM). The dual mesh $\mathcal T^*_h$ satisfies $\overline \Omega  = \bigcup\limits_{k = 1}^{N_s } { \overline{\Omega _k^s }}$, and none of the elements of $\mathcal T^*_h$  overlap.\\
\newline
In order to visualize the dual mesh $\mathcal  T^*_h$, we give examples for several elements $\Omega^s_k \in \mathcal T^*_h$ used in the 2D and 3D bES-FEM and the bFS-FEM. For the 2D and 3D bES-FEM, elements $\Omega^s_k \in \mathcal  T^*_h$ are described in Figure~\ref{fig:smoothingCell_bESFEM}.
Figure~\ref{fig:cell_bESFEM_a} illustrates in 2D an element $\Omega_1^s \in \mathcal  T^*_h$ which has an edge $e_1$ aligned with the domain boundary and an interior element $\Omega_2^s \in \mathcal  T^*_h$ centered along an interior edge $e_2$. Figure~\ref{fig:cell_bESFEM_b} shows an element of the dual mesh in 3D consisting of six tetrahedral elements together with an inner smoothing cell centered along the edge AB. In Figures~\ref{fig:cell_bESFEM_c} and~\ref{fig:cell_bESFEM_d} we give a further 3D example showing a smoothing cell $\Omega^s_k$ associated with edge AB of the boundary $\partial \Omega$.\\
\newline
For bFS-FEM, we also have an example for a smoothing domain $\Omega_k^s \in  \mathcal  T^*_h$. The domain $\Omega_k^s$ associated with the face $k$ is created by simply connecting three nodes B, C, D of the face to the centers H, I of adjacent elements as shown in Figure~\ref{fig:smoothingCell_bFSFEM}.
\begin{figure}[htbp!]
  \centering
  \caption{Two adjacent tetrahedral elements and the smoothing domains (shaded region) formed based on their interface $k$ (BCD) in the 3D FS-FEM.}
  \label{fig:smoothingCell_bFSFEM}
  \end{figure}
{With the dual mesh $\mathcal T^*_h$, the space $\mathcal V_h^{\mathcal B}$  is equipped with the following the inner product, semi-norm and norm (see \cite{GRL}).}\\ 
As a consequence of remarks $3.4$ and $3.5$ in \cite{GRL}, we have relationships between $|\cdot|_{\mathcal V_h^\mathcal B}$ and $|\cdot|_{1}$, and also between $||\cdot||_{\mathcal V_h^\mathcal B}$ and $||\cdot||_{1}$ as follows:
\begin{equation}
\label{phtrinh1a}
|\bold{w}|_{\mathcal V_h^\mathcal B} \le |\bold{w}|_{ 1} \quad \text{and}\quad ||\bold{w}||_{\mathcal V_h^\mathcal B} \le ||\bold{w}||_{ 1} \quad \text{with}~\bold{w}\in \mathcal V_h^\mathcal B \subset [H^1(\Omega)]^d,
\end{equation}
where $H^1(\Omega)$ is a Sobolev space which is endowed with the semi-norm $|.|_{1}$ and norm $||.||_{1}$, defined by the inner product $(\bold{w,v})_1 = \sum\limits_{i = 1}^d {(\bold{w}_i ,\bold{v}_i )_1 }$ (see chapter $3$ in \cite{AF}).\\
\subsection{The third mesh}
\indent Next, a third mesh $\mathcal T^{**}_h$ is constructed by connecting all centroids $\{c_T\}_{T\in\mathcal T_h}$ and midpoints of all edges of $\mathcal T_h$ in 2D, plus  barycenter points of all faces in 3D. The third mesh $\mathcal T^{**}_h$ satisfies $\overline \Omega  = \bigcup\limits_{i = 1}^{N_n } { \overline{V_i}}$, and none of the elements of $\mathcal T^*_h$ overlap. Each element $V_k \in \mathcal T^{**}_h$ is also associated with a vertex  $\bold{x}_k$ of the primal mesh.
\begin{figure}[htbp!]
        \centering
        \begin{subfigure}[b]{6.5cm}
                \caption{a 2D element}
                \label{fig:elementsVk_2D}
        \end{subfigure}%
        \quad
        \begin{subfigure}[b]{6.5cm}
                \caption{a 3D element}
                \label{fig:elementsVk_3D}
        \end{subfigure}
        \caption{Elements $V_k$ of $\mathcal T^{**}_h$.}
        \label{fig:elementsVk}
\end{figure}
Figure~\ref{fig:elementsVk_2D} is an example of an element $V_k \in \mathcal T^{**}$ constructed by connecting centroids $\{c_{T_i}\}_{T_i \in \mathcal T_h}$ and midpoints $\{\bold{x}_{e_i}\}_{i \in \overline{1,6}}$ with edges $\{e_i\}_{i = \overline{1,6}}$ in 2D. Figure~\ref{fig:elementsVk_3D} is another example for an intersecting domain $V_k \cap T$ between $V_k \in \mathcal T^{**}_h$ and $T \in \mathcal T_h$ in 3D. This intersecting domain is made from a set of a vertex $\bold{x}_k$, midpoints $\{\bold{x}_{e_i}\}_{i = \overline{1,3}}$ of edges $\{e_i\}_{i =\overline{1,3}}$, barycentric points  $\{\bold{x}_{f_i}\}_{i = \overline{1,3}}$ of faces $\{ f_i\}_{i=\overline{1,3}}$, a centroid $c_{T}$.\\
Based on this third mesh, we define the following finite element space for the pressure
\begin{displaymath}
\mathcal V^{**}_h = \left\{ {p \in L_0^2 (\Omega )\,\,\text{such that}\,p\left| {_V  \in \mathbb P^0 (V)} \right.,\,V \in \mathcal T_h^{**} \,\,} \right\},
\end{displaymath}
where its norm $||.||_0$ of $\mathcal V^{**}_h $  is  defined by $||q||_0  = \left( {\int\limits_\Omega  {q^2 \text{d}\Omega} } \right)^{\frac{1}{2}}$ for $\forall q \in \mathcal V^{**}_h$.\\
Let $p_i$ be the nodal value of $p_h$ at a vertex $i \in \overline{1,N_n}$.  Then $p_h = \sum\limits_{i = 1}^{N_n } {p_i \chi _i }$ if $p_h \in \mathcal V^{**}_h$, where $\chi_i$ are the characteristic functions of $V_i \in \mathcal T_h^{**}$, $i = \overline{1,N_n}$.\\
\newline
Now, we apply 2D/3D bES-FEM and bFS-FEM for discretizing the nearly-incompressible elasticity problem in the two following sections.
\subsection{Smoothed strain and smoothed divergence}
\label{Smooth1}
In 2D, according to the formula (\ref{phtrinh03}), the discretized strain $\boldsymbol{\varepsilon}(\bold{u}_h)$ is obtained as
\begin{equation}
\label{phtrinh2a}
\boldsymbol{\varepsilon}(\bold{u}_h) = \partial \bold{u}_h =\left[ {\begin{array}{*{20}c}
   {\frac{\partial }{{\partial x}}} & 0  \\
   0 & {\frac{\partial }{{\partial y}}}  \\
   {\frac{\partial }{{\partial y}}} & {\frac{\partial }{{\partial x}}}  \\
\end{array}} \right]\bold{u}_h.
\end{equation}
On each smooth element $\Omega^s_k \in \mathcal T_h^*$, the strain $\boldsymbol{\varepsilon}(\bold{u}_h)$ is smoothed by
\begin{equation}
\label{phtrinh2}
\boldsymbol{\overline \varepsilon}^{(k)}(\bold{u}_h)  = \frac{1}{{m(\Omega _k^s )}}\int\limits_{\Omega _k^s } {\boldsymbol{\varepsilon}(\bold{u}_h (\bold{x}))\text{d}\Omega}  =  \frac{1}{{m(\Omega _k^s )}}\int\limits_{\Omega _k^s } {\partial \bold{u}_h (\bold{x})\text{d}\Omega} \quad\text{with} ~\bold{u}_h \in \mathcal  V^{\mathcal B}_h,
\end{equation}
and we also have a formula for the smoothed divergence
\begin{equation}
\label{phtrinh3}
\left(\overline {\nabla \cdot\bold{u}_h}\right)|_{\Omega^s_k}  = \frac{1}{{m(\Omega _k^s )}}\int\limits_{\Omega _k^s } {\nabla \cdot\bold{u}_h (\bold{x})\text{d}\Omega}  \quad \text{with} ~\bold{u}_h \in \mathcal  V^{\mathcal B}_h.
\end{equation}
By performing the integration in (\ref{phtrinh2}), the smoothed strain $\boldsymbol{\overline \varepsilon}^{k}$ can be rewritten on the boundary $\partial \Omega^k_s$, as follows:
\begin{equation}
\label{phtrinh3a}
\boldsymbol{\overline \varepsilon}^{(k)}(\bold{u}_h)  = \frac{1}{{m(\Omega _k^s )}}\int\limits_{\partial \Omega_k^s } { \bold{n}^{(k)}(\bold{x}) \bold{u}_h (\bold{x}) \text{d}\gamma(\bold{x})}
\end{equation}
where $\bold{n}^{(k)}(\bold{x})$ is defined by $\left[ {\begin{array}{*{20}c}
   {{n}_x^{(k)} } & 0  \\
   0 & {{n}_y^{(k)} }  \\
   {{n}_y^{(k)} } & {{n}_x^{(k)} }  \\
\end{array}} \right]$, and the two notations ${{n}_x^{(k)} }$, ${{n}_y^{(k)} }$ are two elements of the outward normal unit vector on the boundary $\partial \Omega_k^s$.\\
\newline
By transforming (\ref{phtrinh1}), (\ref{phtrinh2a}) and (\ref{phtrinh2}) into the formula (\ref{phtrinh3a}), we remove the need to use shape function derivatives in the calculation of the discrete smoothed strain $\boldsymbol{\overline \varepsilon}^k(\bold{u}_h)$. The number of Gauss points used for the line (2D) or face (3D) integration in (\ref{phtrinh3a}) depends on the order of the shape functions and bubble functions. In 3D, the strain and the divergence are similarly smoothed.\\
\subsection{Weakened weak statement for bES-FEM and bFS-FEM}
\label{smooth2}
Here, we want to find the discrete solution $(\bold{u}_h,p_h) \in \mathcal V^{\mathcal B}_h \times \mathcal V^{**}_h$ such that
\begin{subequations}
\label{phtrinh4}
\begin{align}
\overline a(\bold{u}_h ,\bold{v}_h ) +  \overline b(\bold{v}_h ,p_h )     & =  (\bold{f},\bold{v}_h )  & \forall \bold{v}_h  \in\mathcal V_h^{\mathcal B} ,\\
\overline b(\bold{v}_h ,p_h ) - \frac{1}{\lambda } \overline c(p_h ,q_h ) & = 0                        & \forall q_h  \in \mathcal V_h^{**} , \label{phtrinh4_pressure}
\end{align}
\end{subequations}
{where
\begin{eqnarray*}
\overline a(\bold{u}_h ,\bold{v}_h) & = & 2\mu \sum\limits_{k = 1}^{N_s } m(\Omega _k^s ) (\boldsymbol{\overline\varepsilon} ^{(k)} (\bold{v}_h))^T \, \bold D \, \boldsymbol{\overline\varepsilon}^{(k)} (\bold{u}_h),\quad
\overline b(\bold{v}_h,p_h )       =  \int\limits_\Omega  {(\overline{\nabla \cdot\bold{v}_h})p_h \,\text{d}\Omega},\\
\overline c(p_h ,q_h )              & = & \int\limits_\Omega  {p_h q_h \,\text{d}\Omega}, \quad
 (\bold{f},\bold{v}_h)  =  \int\limits_\Omega  \bold{v}_h^T (\bold{x}) \bold{f}(\bold{x}) \text{d}\Omega.
\end{eqnarray*}}
The system of equations in (\ref{phtrinh4}) is known as a weakened weak ($W^2$) form because derivatives of the displacements are no longer needed in contrast to the usual weak form \cite{LAMI}. Also, due to (\ref{phtrinh4_pressure}), we will be able to calculate the discrete pressure $p_h$ from the smoothed divergence $\overline{\nabla\cdot\bold{u}_h}$ as is shown by the formula in (\ref{phtrinh44a}), see Remark \ref{mathproper}.2.
\section{The mathematical properties}
\label{mathproper}
In this section, we present the important mathematical results for bES-FEM and bFS-FEM when applied to the linear elasticity problem.\\
\newline
{\bf Theorem \ref {mathproper}.1} \textit{(Coercivity and Continuity)}\\
\textit{ The bilinear form $\overline a(\cdot,\cdot)$ is continuous, symmetric and coercive on
\begin{displaymath}
\mathcal V^\mathcal B_{h,0} := \left \{ \bold{v} \in \mathcal V^\mathcal B_h : \overline b(\bold{v},q) = 0, \forall q \in \mathcal V^{**}_h \subset  L^2_0(\Omega) \right \},
\end{displaymath}
i.e. there exists an $\alpha_0, \alpha_1 > 0$ such that}
\begin{equation}
\label{phtrinh5a}
\overline a(\bold{v,v}) \ge \alpha_0 ||\bold{v}||^2_{\mathcal V_h^\mathcal B}, \quad \bold{v} \in \mathcal V^\mathcal B_{h,0}. \quad (\textit{coercivity})
\end{equation}
\begin{equation}
\label{ptrinh5b}
\overline a(\bold{v,w}) \le \alpha_1 ||\bold v||_{\mathcal V^{\mathcal B}_h} ||\bold w||_{\mathcal V^{\mathcal B}_h}, \quad \bold v, \bold w \in \mathcal V^{\mathcal B}_{h,0} \quad (\textit{continuity})
\end{equation}
This theorem can be proven by invoking the theorem $3.2$ (\textit{coercivity}) and the theorem $3.3$ (\textit{continuity}) in \cite{GRL1}. $\hspace{12.5cm} \square$\\
\newline
{\bf Theorem \ref {mathproper}.2} \textit{(Stability)}\\
\textit{The bilinear form $\overline b(\cdot,\cdot)$ on $\mathcal V_h^{\mathcal B} \times \mathcal V_h^{**}$ is continuous and satisfies the uniform inf-sup condition, i.e. there exists a positive constant $\beta_0$ independent of the mesh size such that
\begin{equation}
\label{phtrinh5}
\mathop {\sup }\limits_{\bold{u}_h  \in \mathcal V_h^ \mathcal B, \bold u_h \ne \bold 0 } \frac{{\overline b(\bold{u}_h,q_h )}}{{||\bold{u}_h ||_{\mathcal V_h^\mathcal B} }} {\quad \mathop {\ge}\limits_{\text{because of} ~(\ref{phtrinh1a})} \quad} \mathop {\sup }\limits_{\bold{u}_h  \in \mathcal V_h^ \mathcal B, \bold u_h \ne \bold 0 } \frac{{\overline b(\bold{u}_h ,q_h )}}{{||\bold{u}_h ||_{1} }} \ge \beta_0 ||q_h ||_0 ,\quad \quad q_h  \in \mathcal  V_h^{**}.
\end{equation}}
To prove the theorem \ref{mathproper}.2, we need to look for a relationship between  $\overline b(\bold{u}_h ,q_h )$ and $b(\bold{u}_h ,q_h ) = \int\limits_\Omega  {\nabla \cdot\bold{u}_h (\bold{x})q_h (\bold{x})\text{d}\Omega}$ with $(\bold{u}_h, q_h) \in  \mathcal V^\mathcal B_h \times \mathcal V^{**}_h$. In \cite{LAMI}, $b(\bold{u}_h ,q_h )$ satisfies the uniform inf-sup condition, from which it follows that $\overline b(\bold{u}_h ,q_h )$ satisfies this condition. This idea was similarly used to prove the uniform inf-sup condition in \cite{LAMI}, where the author also indicated the relationship between $b(\bold{u}_h ,q_h )$ and the bilinear form derived for the MINI element.\\
\newline
Let $(\bold{u}_h,q_h) \in \mathcal V_h^\mathcal B \times \mathcal V^{**}_h$, we have
\begin{equation}
\label{phtrinh6}
\overline b(\bold{u}_h ,q_h ) = \int\limits_\Omega  {(\overline{\nabla \cdot\bold{u}_h})q_h (\bold{x})\text{d}\Omega} = \int\limits_\Omega  {(\overline {\nabla \cdot\boldsymbol{\ell}_h}  + \overline{\nabla \cdot\bold{b}_h} )\,q_h (\bold{x})\text{d}\Omega},
\end{equation}
where there exists uniquely $\boldsymbol{\ell}_h \in \mathcal V_h$ and $\bold{b}_h \in [\mathcal B_h]^d$ such that $\bold{u}_h = \boldsymbol{\ell}_h + \bold{b}_h$. In (\ref{phtrinh6}), the smoothed divergences  $\overline {\nabla \cdot\boldsymbol{\ell}_h}$ and $\overline{\nabla \cdot\bold{b}_h} $, which are restricted on  $\Omega_k^s \in \mathcal T^*_h$, are defined  by (\ref{phtrinh3}).\\

\noindent {\bf Lemma \ref{mathproper}.1} \textit{The value of }$\int\limits_\Omega  {(\overline {\nabla \cdot\boldsymbol{\ell} _h })q_h (\bold{x}) \text{d}\Omega}  - \int\limits_\Omega  {(\nabla \cdot\boldsymbol{\ell} _h )q_h (\bold{x}) \text{d}\Omega}$ \textit{is equal to $0$ with $(\boldsymbol{\ell}_h, q_h) \in \mathcal V_h \times \mathcal V^{**}_h$}.\\
\newline
\textit{Proof:} \quad Using the fact that $\nabla \cdot\boldsymbol{\ell} _h$ is constant on each $T \in \mathcal T_h$, we obtain
\begin{equation}
\label{phtrinh7}
\int\limits_\Omega  {(\nabla \cdot\boldsymbol{\ell} _h) \,q_h (\bold{x})\text{d}\Omega}  = \sum\limits_{T \in \mathcal T_h } {(\nabla \cdot\boldsymbol{\ell} _h )\left| {_T } \right.\int\limits_T {q_h (\bold{x})\text{d}\Omega} }.
\end{equation}
For any element $T \in \mathcal T_h$ with its vertices $\{\bold{x}_{T^{(i)}}\}_{i=\overline{1,d+1}}$, we have
\begin{equation}
\label{phtrinh8}
(\nabla \cdot\boldsymbol{\ell}_h )\left| {_T } \right.\int\limits_T {q_h (\bold{x})\text{d}\Omega}  = (\nabla \cdot\boldsymbol{\ell}_h )\sum\limits_{i = 1}^{d + 1} {m\left( {V_{\bold{x}_{T^{(i)} } }  \cap T} \right) q_{T^{(i)} } } =  (\nabla \cdot\boldsymbol{\ell}_h )\left| {_T } \right.\sum\limits_{i = 1}^{d + 1} {\frac{{m(T)}}{d+1} q_{T^{(i)} } },
\end{equation}
where for each $i = \overline{1,d+1}$, $V_{\bold{x}_{T^{(i)}}} \in \mathcal T^{**}_h$ is  associated with a vertex  $\bold{x}_{T^{(i)}}$ of $T$, and ${q_{T^{(i)} } }$ is a nodal value of $q_h$ at  a vertex ${\bold{x}_{T^{(i)}}}$. The value ${m({V_{\bold{x}_{T^{(i)} } }  \cap T)}}$ is equal to $\frac{{m(T)}}{d+1}$ with $i =\overline{1,d+1}$, because  the third mesh $\mathcal T_h^{**}$ is constructed by barycentric points of all faces (3D), midpoints of all edges and the centroid points $c_T$ for all $T \in \mathcal T_h$. \\
\newline
We now calculate the integral $\int\limits_\Omega  {(\overline{\nabla \cdot\boldsymbol{\ell}_h}) \,q_h (\bold{x})\text{d}\Omega}$, for the two methods, bES-FEM and bFS-FEM, as follows:\\
\newline
\noindent {\bf For the 2D and 3D bES-FEM}\\
\newline
In (\ref{phtrinh6}), we  consider
\begin{equation}
\label{phtrinh9}
\int\limits_\Omega  {(\overline{\nabla \cdot\boldsymbol{\ell}_h}) \,q_h (\bold{x})\text{d}\Omega}  = \sum\limits_{T \in \mathcal T_h } {\int\limits_T {(\overline{\nabla \cdot\boldsymbol{\ell} _h} ) q_h (\bold{x})\text{d}\Omega} }.
\end{equation}
On the above element $T \in \mathcal T_h$, the integral $\int\limits_T {(\overline{\nabla \cdot\boldsymbol{\ell}_h} )q_h (\bold{x})\text{d}\Omega}$ is computed by
\begin{equation}
\label{phtrinh10}
\int\limits_T {(\overline{\nabla \cdot\boldsymbol{\ell}_h} ) q_h (\bold{x})\text{d}\Omega} =\sum\limits_{i = 1}^{d + 1} {\left[ {\sum\limits_{\scriptstyle \,\,\,\,\,\,\,\,\, \hfill \atop
  \scriptstyle e_{T^{(i)} }  \in \mathcal E_{T^{(i)} }  \hfill} {\frac{{m\left( {V_{\bold{x}_{T^{(i)} } }  \cap T \cap \Omega _{e_{T^{(i)} } }^s } \right)}}{{m\left( {\Omega _{e_{T^{(i)} }}^s } \right)}}} \int\limits_{\Omega _{e_{T^{(i)} } }^s } {\nabla \cdot\boldsymbol{\ell}_h \text{d}\Omega} } \right]} q_{T^{(i)} },
\end{equation}
where the domain ${\Omega _{e_{T^{(i)} }}^s } \in \mathcal T^*_h$ corresponds to the edge $e_{T^{(i)} }$. The set  $\mathcal E_{{T^{(i)}}}$ contains all edges of $T$ such that these edges have a common vertex $\bold{x}_{T^{(i)}}$.  \\
\newline
In the first case of $T$ (a triangle or tetrahedron), we assume that all edges and all faces (3D) of $T$ are inner edges and inner faces, {\it i.e.} the edges and faces are not on the boundary $\partial \Omega$. For each $i = \overline{1,d+1}$ and $j = \overline{1,d}$, the integral $\int\limits_{\Omega _{e_{T^{(i)} }}^s } {\nabla \cdot\boldsymbol{\ell}_h \text{d}\Omega} $  is computed by
\begin{equation}
\label{phtrinh11}
\int\limits_{\Omega _{e_{T^{(i)} } }^s } {\nabla \cdot\boldsymbol{\ell}_h \text{d}\Omega}  = m\left( {\Omega _{e_{T^{(i)} } }^s  \cap T} \right) \left( {\nabla \cdot\boldsymbol{\ell}_h } \right)\left| {_T } \right. + \sum\limits_{K \in \mathcal T_{e_{T^{(i)} }} \backslash \{ T\} } {m\left( {\Omega _{e_{T^{(i)} } }^s  \cap K} \right) \left( {\nabla \cdot\boldsymbol{\ell}_h } \right)\left| {_K } \right.} ,
\end{equation}
where $\mathcal T_{e_{T^{(i)} }}$ is a subset of $\mathcal T_h$ such that its elements have a common edge $e_{T^{(i)} }$ and $T \in \mathcal T_{e_{T^{(i)} }}$. \\
From  (\ref{phtrinh10}) and (\ref{phtrinh11}), the integral $\int\limits_{T}{(\overline{\nabla \cdot\boldsymbol{\ell} _h} ) q_h (\bold{x})\text{d}\Omega}$  has the coefficient of $(\nabla \cdot\boldsymbol{\ell}_h )\left| {_T} \right. q_{T^{(i)}}$
\begin{equation}
\label{phtrinh13}
\sum\limits_{\scriptstyle \,\,\,\,\,\,\,\,\, \hfill \atop
  \scriptstyle e_{T^{(i)} }  \in \mathcal E_{T^{(i)} }  \hfill} {\frac{{m\left( {V_{\bold{x}_{T^{(i)} } }  \cap T \cap \Omega _{e_{T^{(i)} } }^s } \right)\,m\left( {\Omega _{e_{T^{(i)} } }^s  \cap T} \right)}}{{m\left( {\Omega _{e_{T^{(i)} } }^s } \right)}}}.
\end{equation}
Together $\int\limits_{T}{(\overline{\nabla \cdot\boldsymbol{\ell}_h} ) q_h (\bold{x})\text{d}\Omega}$ , we only find the coefficient of $(\nabla \cdot\boldsymbol{\ell}_h )\left| {_T} \right. q_{T^{(i)}}$ in  $\int\limits_{K \in \mathcal T_{e_{T^{(i)} } } \backslash \{ T\} }{(\overline{\nabla \cdot\boldsymbol{\ell}_h} ) q_h (\bold{x})\text{d}\Omega}$ for all $K \in \mathcal T_{e_{T^{(i)} } } \backslash \{ T\}$ and $e_{T^{(i)} }  \in \mathcal E_{T^{(i)} }$, as follows:
\begin{equation}
\label{phtrinh14}
\left\{ {\frac{{m\left( {V_{\bold{x}_{T^{(i)} } }  \cap K \cap \Omega _{e_{T^{(i)} } }^s } \right)\,m\left( {\Omega _{e_{T^{(i)} } }^s  \cap T} \right)}}{{m\left( {\Omega _{e_{T^{(i)} } }^s } \right)}}} \right\}_{\forall K \in \mathcal T_{e_{T^{(i)} } }\backslash \{ T\} ~\text{and}~\forall e_{T^{(i)} }  \in \mathcal E_{T^{(i)} }}.
\end{equation}
From (\ref{phtrinh13}) and (\ref{phtrinh14}), in the integral $\int\limits_\Omega  {(\overline{\nabla \cdot\boldsymbol{\ell}_h} )q_h (\bold{x})\text{d}\Omega}$, the coefficient of  $(\nabla \cdot\boldsymbol{\ell}_h )\left| {_T} \right. q_{T^{(i)}}$ is equal to {
\begin{equation}
\label{phtrinh15}
\sum\limits_{\,\,\,\,\,\,\,\,\,e_{T^{(i)} }  \in \mathcal E_{T^{(i)} } } {\left[ \begin{array}{l}
 \frac{{m\left( {V_{\bold{x}_{T^{(i)} } }  \cap T \cap \Omega _{e_{T^{(i)} } }^s } \right)\,m\left( {\Omega _{e_{T^{(i)} } }^s  \cap T} \right)}}{{m\left( {\Omega _{e_{T^{(i)} } }^s } \right)}} +
 \sum\limits_{K \in \mathcal T_{e_{T^{(i)} } } \backslash \{ T\} } {\frac{{m\left( {V_{\bold{x}_{T^{(i)} } }  \cap K \cap \Omega _{e_{T^{(i)} } }^s } \right)\,m\left( {\Omega _{e_{T^{(i)} } }^s  \cap T} \right)}}{{m\left( {\Omega _{e_{T^{(i)} } }^s } \right)}}}  \\
 \end{array} \right]}.
\end{equation}}
By using the centroids $c_T$ for all $T \in \mathcal T_h$, the midpoints of all edges, plus  barycentric points of all faces (3D) to construct the dual mesh $\mathcal T^*_h$ and the third mesh $\mathcal T^{**}_h$, we have
{\begin{eqnarray}\label{phtrinh15a}
m\left( {V_{\bold{x}_{T^{(i)} } }  \cap T \cap \Omega _{e_{T^{(i)} } }^s } \right) & = & \frac{{m(T)}}{{\card(\mathcal E_{{T^{(i)} }})(d + 1)}}=\frac{{m(T)}}{{d(d + 1)}},  ~
m\left( {\Omega _{e_{T^{(i)} } }^s  \cap T} \right)  =  \frac{{m(T)}}{{\card(\mathcal E_T)}},  \nonumber \\
m\left( {V_{\bold{x}_{T^{(i)} } }  \cap K \cap \Omega _{e_{T^{(i)} } }^s } \right)  &=& \frac{{m(K)}}{{\card(\mathcal E_{{K^{(i)} }})(d + 1)}} = \frac{{m(K)}}{{d(d + 1)}}, ~
m\left( {\Omega _{e_{T^{(i)} } }^s } \right)    =  \sum\limits_{L \in \mathcal T_{e_{T^{(i)} } } } {\frac{{m(L)}}{{\card(\mathcal E_L)}}}, \nonumber \\
\end{eqnarray}}
where  for all $i =\overline{1,d}$ and $T \in \mathcal T_h$, the notations $\card(\mathcal E_{T^{(i)} })$ and $\card(\mathcal E_T)$  are the number of all elements of $\mathcal E_{T^{(i)} }$ and $\mathcal E_T$, respectively. Furthermore, we have $\card(\mathcal E_{T^{(i)} }) = d$ and $\card(\mathcal E_T) = \card(\mathcal E_K)$ for all $K, T \in \mathcal T_h$, because the primal mesh $\mathcal T_h$ is a triangulation\\
\newline
Therefore, the coefficient of  $(\nabla \cdot v_h )\left| {_T} \right. q_{T^{(i)}}$ is computed as
\begin{eqnarray}
\label{phtrinh16}
\sum\limits_{\,\,\,\,\,\,\,\,\,e_{T^{(i)} }  \in \mathcal E_{T^{(i)} } } {\left[ {\frac{{\frac{{m(T)}}{{\card(\mathcal E_{{T^{(i)} }})(d + 1)}}\,\frac{m(T)}{\card(\mathcal E_T)}}}{{\sum\limits_{L \in \mathcal T_{e_{T^{(i)} } } } {\frac{{m(L)}}{{\card(\mathcal E_L)}}} }} + \sum\limits_{K \in \mathcal T_{e_{T^{(i)} } } \backslash \{ T\} } {\frac{{\frac{{m(K)}}{{\card(\mathcal E_{{K^{(i)} }})(d + 1)}}\,\frac{m(T)}{\card(\mathcal E_K)}}}{{\sum\limits_{L \in \mathcal T_{e_{T^{(i)} } } } {\frac{{m(L)}}{{\card(\mathcal E_L)}}} }}} } \right]} = \frac{{m(T)}}{{d + 1}}.
\end{eqnarray}
From (\ref{phtrinh8}) and (\ref{phtrinh16}), the two coefficients of  $(\nabla \cdot\boldsymbol{\ell}_h )\left| {_T} \right. q_{T^{(i)}}$  in the two integrals $\int\limits_\Omega  {(\overline{\nabla \cdot\boldsymbol{\ell}_h} )q_h (\bold{x})\text{d}\Omega}$ and $\int\limits_\Omega  {({\nabla \cdot\boldsymbol{\ell}_h} )q_h (\bold{x})\text{d}\Omega}$ are equal.\\
\newline
\noindent {\bf For the bFS-FEM method}\\
\newline
Using this method, we obtain the coefficient of $(\nabla\cdot\boldsymbol{\ell}_h)|_T q_{T^{(i)}}$
in the integral $\int\limits_\Omega  {\overline{(\nabla \cdot\boldsymbol{\ell}_h )}q_h (\bold{x})\text{d}\Omega}$ to be
{\begin{equation}
\label{phtrinh16a}
\sum\limits_{\scriptstyle   f_{T^{(i)} }  \in \mathcal  F_{T^{(i)} } , ~\Omega _{f_{T^{(i)} } }^s  \in \mathcal T^{**}_h \hfill \atop  \scriptstyle K \in \mathcal T_h ,~\mathcal F_K  \cap \mathcal F_T  = \{ f_{T^{(i)} } \} \hfill} {\left[ \begin{array}{l}
 \frac{{m\left( {V_{\bold{x}_{T^{(i)} } }  \cap T \cap \Omega _{f_{T^{(i)} } }^s } \right)\,m\left( {\Omega _{f_{T^{(i)} } }^s  \cap T} \right)}}{{m\left( {\Omega _{f_{T^{(i)} } }^s } \right)}} +
 \frac{{m\left( {V_{\bold{x}_{T^{(i)} } }  \cap K \cap \Omega _{f_{T^{(i)} } }^s } \right)\,m\left( {\Omega _{f_{T^{(i)} } }^s  \cap T} \right)}}{{m\left( {\Omega _{f_{T^{(i)} } }^s } \right)}} \\
 \end{array} \right]} =\frac{m(T)}{d+1},
\end{equation}}
where $\mathcal F_{T^{(i)}}$ is a set of all faces of a tetrahedral $T$ whose has a common vertex $\bold{x}_{T^{(i)}}$ and $\card(\mathcal F _{T^{(i)}})$ is equal to $d$. The notation $f_{T^{(i)}}$ is a face of $T$, one of its vertices is $\bold{x}_{T^{(i)}}$. The two sets $\mathcal F_K$, $\mathcal F_T$ contain all faces of $K, T \in \mathcal T_h$, respectively. We have used the following expressions
\begin{align*}
\label{phtrinh16b}
m\left( {V_{\bold{x}_{T^{(i)} } }  \cap T \cap \Omega _{f_{T^{(i)} } }^s } \right) &= \frac{{m(T)}}{{d(d + 1)}}, & m\left( {\Omega _{f_{T^{(i)} } }^s  \cap T} \right) & = \frac{{m(T)}}{{d + 1}},  \\[2mm]
m\left( {V_{\bold{x}_{T^{(i)} } }  \cap K \cap \Omega _{f_{T^{(i)} } }^s } \right) &= \frac{{m(K)}}{{d(d + 1)}}, & m\left( {\Omega _{f_{T^{(i)} } }^s } \right)        & = \frac{{m(T) + m(K)}}{{d + 1}}.
\end{align*}
In the other cases of  $T \in \mathcal T_h$ which has at least one edge or one face belonging to the boundary $\partial \Omega$, we also obtain the same results as (\ref{phtrinh16}) and (\ref{phtrinh16a}).\\
\newline
From (\ref{phtrinh8}), (\ref{phtrinh16}) and (\ref{phtrinh16a}), we deduce that
\begin{displaymath}
\int\limits_\Omega  {(\overline {\nabla \cdot\boldsymbol{\ell}_h })q_h (\bold{x})\text{d}\Omega}  - \int\limits_\Omega  {(\nabla \cdot\boldsymbol{\ell}_h )q_h (\bold{x})\text{d}\Omega} = 0. \hspace{4cm} \square
\end{displaymath}
{\bf Remark 4.1:} Due to the result from \cite{LAMI} and Lemma \ref{mathproper}.1,
we can conclude that if the displacement space is not enriched by bubble functions, the 2D/3D ES-FEM and the 3D FS-FEM  violate the uniform inf-sup condition, further discussed in Remark $4.2$.\\
\newline
Our next objective is to find the relationship between $\int\limits_\Omega  {(\overline{\nabla \cdot\bold{b}_h} )q_h (x)\text{d}\Omega}$ and $\int\limits_\Omega  {(\nabla \cdot\bold{b}_h )q_h (\bold{x})\text{d}\Omega}$. This relationship is shown in the following lemma.\\
\newline
{\bf Lemma \ref{mathproper}.2} \textit{There exists a positive constant $\alpha$ which depends on the bubble function, such that }
\begin{displaymath}
\int\limits_\Omega  {\overline{\nabla \cdot\bold{b}_h (\bold{x})}q_h (\bold{x}) \text{d}\Omega} = \alpha \int\limits_\Omega  {\nabla \cdot\bold{b}_h (\bold{x})q_h (\bold{x})\text{d}\Omega}.
\end{displaymath}
\textit{Proof:} \quad By the definitions of the spaces $\mathcal B_h$ and $\mathcal V^{**}_h$, with $(\bold{b}_h,q_h) \in \mathcal B_h \times \mathcal V^{**}_h$,  we get
\begin{equation}
\label{phtrinh17}
\int\limits_\Omega  {\nabla \cdot\bold{b}_h (\bold{x})q_h (\bold{x})\text{d}\Omega}  = \sum\limits_{T \in \mathcal  T_h } \quad{\sum\limits_{\scriptstyle i = 1,\,V_i  \in \mathcal  T_h^{**}  \hfill \atop
  \scriptstyle  T \cap V_i  \ne \emptyset  \hfill}^{N_n } \left[ {\ \int\limits_{V_i  \cap T} {\nabla \cdot\bold{b}_h (\bold{x})q_i \text{d}\Omega} } \right]}
\end{equation}
and
\begin{equation}
\label{phtrinh18}
\int\limits_\Omega  {\overline{\nabla \cdot\bold{b}_h (\bold{x})}q_h (\bold{x})\text{d}\Omega}  = \sum\limits_{T \in \mathcal T_h }\quad {\sum\limits_{\scriptstyle   i = 1,\,V_i  \in \mathcal  T_h^{**}  \hfill \atop
  \scriptstyle  T \cap V_i  \ne \emptyset  \hfill}^{N_n }  \left[\ {\int\limits_{V_i  \cap T} {\overline{\nabla \cdot\bold{b}_h (\bold{x})}q_i \text{d}\Omega} } \right]}.
\end{equation}
Considering $T$ whose all edges stay in the internal domain $\Omega$, for each $i = \overline{1,d+1}$, we have
\begin{equation}
\label{phtrinh19}
\int\limits_{V_{\bold{x}_{T^{(i)}} }  \cap T} {\nabla \cdot\bold{b}_h (\bold{x})q_{T^{(i)}} \text{d}\Omega}  = q_{T^{(i)}} \bold{u}_{c_T } \cdot \hspace{-0.3cm}\int\limits_{V_{\bold{x}_{T^{(i)}} }  \cap T} {\nabla N_{c_T }^b (\bold{x})\text{d}\Omega},
\end{equation}
where $\bold{b}_h$ is rewritten as $\bold{b}_h = \bold{u}_{c_T} N_{c_T }^b $. Once again the calculation of the integral $\int\limits_\Omega  {\overline{\nabla \cdot\bold{b}_h (\bold{x})}q_h (\bold{x})\text{d}\Omega}$ is performed for bES-FEM and bFS-FEM in turn.\\
\newline
{\bf For the 2D and 3D bES-FEM}
{\begin{equation}
\label{phtrinh20}
\begin{split}
&\int\limits_{V_{\bold{x}_{T^{(i)}}}  \cap T} {\overline{\nabla \cdot\bold{b}_h (\bold{x})}q_{T^{(i)}} \text{d}\Omega}  = \left[ {\sum\limits_{\,\,\,\,\,\,\,\,\,e_{T^{(i)} } \in \mathcal E_{T^{(i)} } } {m\left( {V_{\bold{x}_{T^{(i)} } }  \cap T \cap \Omega _{e_{T^{(i)} } }^s } \right)\left. {\left( {\nabla \cdot\bold{b}_h } \right)} \right|_{\Omega _{e_{T^{(i)} } }^s } } } \right]q_{T^{(i)} }  \\
&= \left\{ {\sum\limits_{\,\,\,\,\,\,\,\,\,e_{T^{(i)} }  \in \mathcal  E_{T^{(i)} } } {\frac{{m\left( {V_{\bold{x}_{T^{(i)} } }  \cap T \cap \Omega _{e_{T^{(i)} } }^s } \right)}}{{m\left( {\Omega _{e_{T^{(i)} } }^s } \right)}}\,\left[ {\sum\limits_{K \in \mathcal T_{e_{T^{(i)} } }  \subset \, \mathcal T_h} {\bold{u}_{c_K } \cdot\left( {\int\limits_{\Omega _{e_{T^{(i)} } }^s  \cap K} {\nabla N_{c_K}^b (\bold{x})\text{d}\Omega} } \right)} } \right]} } \right\}q_{T^{(i)} }.
\end{split}
\end{equation}}
From (\ref{phtrinh20}),  with $T \in  \mathcal T_{e_{T^{(i)}}}$, the coefficient of $q_{T^{(i)}} \bold{u}_{c_T}$  in $\int\limits_{V_{\bold{x}_{T^{(i)}} }  \cap T} {(\overline{\nabla \cdot\bold{b}_h})q_{T^{(i)}} \text{d}\Omega}$ is equal to
\begin{equation}
\label{phtrinh21}
\sum\limits_{\,\,\,\,\,\,\,\,\,e_{T^{(i)} }  \in \mathcal E_{T^{(i)} } } {\left[ {\frac{{m\left( {V_{\bold{x}_{T^{(i)} } }  \cap T \cap \Omega _{e_{T^{(i)} } }^s } \right)}}{{m\left( {\Omega _{e_{T^{(i)} } }^s } \right)}}\int\limits_{\Omega _{e_{T^{(i)} } }^s  \cap \,T} {\nabla N_{c_T }^b (\bold{x})\text{d}\Omega} } \right]}.
\end{equation}
Furthermore, the other coefficients of $q_{T^{(i)}} \bold{u}_{c_T}$, which are also found in\\
$ 
\left\{ {\int\limits_{V_{\bold{x}_{T^{(i)}} }  \cap K} {(\nabla \cdot\bold{b}_h )q_{T^{(i)} } \text{d}\Omega} } \right\}_{K \in  \mathcal T_{e_{T^{(i)} } } \backslash  \{T \} \subset \, \mathcal T_h, \,\forall e_{T^{(i)} }  \in \mathcal E_{T^{(i)} } }~\text{are equal to}
$
\begin{equation}
\label{phtrinh22}
\left\{ {\frac{{m\left( {V_{\bold{x}_{T^{(i)} } }  \cap K \cap \Omega _{e_{T^{(i)} } }^s } \right)}}{{m\left( {\Omega _{e_{T^{(i)} } }^s } \right)}}\int\limits_{\Omega _{e_{T^{(i)} } }^s  \cap \,T} {\nabla N_{c_T }^b (\bold{x})\text{d}\Omega} } \right\}_{\scriptstyle K \in  \mathcal T_{e_{T^{(i)} } } \backslash  \{T\} \subset \mathcal T_h  \hfill \atop
  \scriptstyle \forall e_{T^{(i)} }  \in \mathcal  E_{T^{(i)} } \, \hfill}.
\end{equation}
From equations (\ref{phtrinh15a}), (\ref{phtrinh21}) and (\ref{phtrinh22}), the coefficient of  $q_{{T^{(i)} }} \bold{u}_{c_T}$ in $\int\limits_\Omega  {(\overline{\nabla \cdot\bold{b}_h})q_h (\bold{x})\text{d}\Omega}$ is given by
{\begin{eqnarray}
\label{phtrinh23}
&&\sum\limits_{\,\,\,\,\,\,\,\,\,e_{T^{(i)} }  \in \mathcal E_{T^{(i)} } } {\left \{ \begin{array}{l}
 \frac{{m\left( {V_{\bold{x}_{T^{(i)} } }  \cap T \cap \Omega _{e_{T^{(i)} } }^s } \right)}}{{m\left( {\Omega _{e_{T^{(i)} } }^s } \right)}} +
 \sum\limits_{K \in  {\mathcal T_{e_{T^{(i)} } } \backslash \{ T}\} \subset \mathcal T_h } {\left [  {\frac{{m\left( {V_{\bold{x}_{T^{(i)} } }  \cap K \cap \Omega _{e_{T^{(i)} } }^s } \right)}}{{m\left( {\Omega _{e_{T^{(i)} } }^s } \right)}}} \right ] }\\
\\
 \end{array} \right\}\int\limits_{\Omega _{e_{T^{(i)} } }^s  \cap \,T} {\nabla N_{c_T }^b (\bold{x})\text{d}\Omega} } \nonumber \\
 \nonumber\\
&& = \frac{1}{d}\sum\limits_{\,\,\,\,\,\,\,\,\,e_{T^{(i)} }  \in \mathcal E_{T^{(i)} } } {\left( {\int\limits_{\Omega _{e_{T^{(i)} } }^s  \cap \,T} {\nabla N_{c_T }^b (\bold{x})\text{d}\Omega} } \right)}.
\end{eqnarray}}

In two dimensions, we compute the coefficient of  $q_{{T^{(i)} }} \bold{u}_{c_T}$ (\ref{phtrinh23}) on the following triangle $T$ having three vertices $\{\bold{x}_i,\bold{x}_j,\bold{x}_k\}$.
\begin{figure}[htbp!]
        \centering
        \begin{subfigure}[b]{6cm}
                \caption{$T_m \cap \Omega^s_{[\bold{x}_k,\bold{x}_l]}$, $T_n \cap \Omega^s_{[\bold{x}_i,\bold{x}_j]}$, $T_l \cap \Omega^s_{[\bold{x}_k,\bold{x}_j]}$}
                \label{fig:intersectingDomains_1}
        \end{subfigure}%
        \quad
        \begin{subfigure}[b]{5cm}
                \caption{$V_{\bold{x}_i}\cap T \cap \Omega^s_{[\bold{x}_i,\bold{x}_j]}$}
                \label{fig:intersectingDomains_2}
        \end{subfigure}
        \caption{Intersecting domains, where $T_m$, $T_n$, $T_l \in \mathcal T_h$;\quad  $\Omega^s_{[\bold{x}_i,\bold{x}_j]}$, $\Omega^s_{[\bold{x}_k,\bold{x}_j]}$,  $\Omega^s_{[\bold{x}_k,\bold{x}_l]} \in \mathcal T^*_h$; and $V_{\bold{x}_i} \in \mathcal T_h^{**}$.}
        \label{fig:intersectingDomains}
\end{figure}
This coefficient is equal to
\begin{eqnarray}
\label{phtrinh23b}
&&\frac{1}{2}\left( {\int\limits_{\Omega _{[\bold{x}_i ,\bold{x}_j ]}^s  \cap T} {\nabla    N_{c_T }^b (\bold{x})\text{d}\Omega}  + \int\limits_{\Omega _{[\bold{x}_k ,\bold{x}_i ]}^s  \cap T} {\nabla    N_{c_T }^b (\bold{x})\text{d}\Omega} } \right) = \nonumber \\
&& \frac{1}{2} \left({\int\limits_{\gamma _j^{(1)} } { N_{c_T }^b (\bold{x})\,\bold{n}_{\gamma _j^{(1)} } \text{d}\gamma (\bold{x})}  + \int\limits_{\gamma _k^{(1)} } { N_{c_T }^b (\bold{x})\,\bold{n}_{\gamma _k^{(1)} } \text{d}\gamma (\bold{x})}} \right).
\end{eqnarray}
\begin{figure}[htbp!]
 \centering
 \caption{A triangular element $(\bold{x}_k,\bold{x}_i,\bold{x}_j)$ of the primal mesh $\mathcal T_h$.}
 \label{fig:elementPrimal}
\end{figure}
In Figure~\ref{fig:elementPrimal}, we introduce some extra notation including the midpoints of edges $[\bold{x}_i,\bold{x}_j]$, $[\bold{x}_k,\bold{x}_i]$ and $[\bold{x}_k,\bold{x}_j]$ denoted by $\bold{x}_{ij}$, $\bold{x}_{ki}$ and $\bold{x}_{kj}$ respectively. We write  $\gamma^{(1)}_k$, $\gamma^{(2)}_k$, $\gamma^{(1)}_j$ and $\gamma^{(2)}_j$ to represent the edges $[\bold{x}_k,c_T]$,  $[\bold{x}_{ij},c_T]$, $[\bold{x}_j,c_T]$ and  $[\bold{x}_{ki},c_T]$. Vectors $\bold{n}_{\gamma^{(1)}_k}$, $\bold{n}_{\gamma^{(2)}_k}$, $\bold{n}_{\gamma^{(1)}_j}$ and $\bold{n}_{\gamma^{(2)}_j}$ are outward normal vectors of $\Omega^s_{[\bold{x}_i,\bold{x}_k]} \cap T$, $V_{\bold{x}_i} \cap T$, $\Omega^s_{[\bold{x}_i,\bold{x}_j]} \cap T$ and $V_{\bold{x}_i} \cap T$ respectively. The length of each vector $\bold{n}_{\gamma^{(1)}_k}$, $\bold{n}_{\gamma^{(2)}_k}$, $\bold{n}_{\gamma^{(1)}_j}$ and $\bold{n}_{\gamma^{(2)}_j}$ is equal to the length of each segment  $\gamma^{(1)}_k$, $\gamma^{(2)}_k$, $\gamma^{(1)}_j$ and $\gamma^{(
2)}_j$, so $\bold{n}_{\gamma^{(1)}_k} = 2 \bold{n}_{\gamma^{(2)}_k}$ and $\bold{n}_{\gamma^{(1)}_j} = 2 \bold{n}_{\gamma^{(2)}_j}$, because \
the length of segments $\gamma_k^{(1)}$ and $\gamma_j^{(1)}$ is equal to twice the length of $\gamma_k^{(2)}$ and $\gamma_j^{(2)}$ respectively.\\
\newline
We directly compute the coefficient (\ref{phtrinh23b}) for the two types of bubble functions investigated here.
\begin{itemize}
\renewcommand{\labelitemi}{$\blacksquare$}
\item \textit{The $\xi$th-power bubble functions (\ref{cubic}) with $\xi = 3$ (the cubic bubble functions)}\\
\end{itemize}
Assume that $ \hat T$ is the reference triangle, $M_T$ is the Jacobian of transformation from the triangle $T$ to $\hat T$, $J_T = \det(M_T)$,
{\begin{align*}
\hat \theta_2^{(1)}  = \int\limits_{\hat \gamma_2^{(1)} } { N_{c_{\hat T} }^b (\bold{x})~\text{d}\gamma (\bold{x})},  \hat \theta_2^{(2)}  = \int\limits_{\hat \gamma_2^{(2)} } { N_{c_{\hat T} }^b (\bold{x})\,\text{d}\gamma (\bold{x})},~
\hat \theta_3^{(1)}  = \int\limits_{\hat \gamma_3^{(1)} } { N_{c_{\hat T} }^b (\bold{x})\,\text{d}\gamma (\bold{x})},~ \hat \theta_3^{(2)}  = \int\limits_{\hat \gamma_3^{(2)} } { N_{c_{\hat T} }^b (\bold{x})\,\text{d}\gamma (\bold{x})}
\end{align*}}
 where the notation  $\hat \gamma_i^{(j)}$ represents
\begin{displaymath}
\hat \gamma_1^{(1)} = [\bold{x}_{\hat T^{(1)}}, c_{\hat T}], ~\hat \gamma_1^{(2)} = [\bold{x}_{\hat T^{(23)}}, c_{\hat T}],~
\hat \gamma_2^{(1)} = [\bold{x}_{\hat T^{(2)}}, c_{\hat T}],
\end{displaymath}
\begin{displaymath}
\hat \gamma_2^{(2)} = [\bold{x}_{\hat T^{(13)}}, c_{\hat T}] ,~\hat \gamma_3^{(1)} = [\bold{x}_{\hat T^{(3)}}, c_{\hat T}], ~ \hat \gamma_3^{(2)} = [\bold{x}_{\hat T^{(12)}}, c_{\hat T}]
\end{displaymath}
with points $\bold{x}_{\hat T^{(1)}}  (0,1)$, $\bold{x}_{\hat T^{(2)}}  (0,0)$, $\bold{x}_{\hat T^{(3)}}  (1,0)$, $\bold{x}_{\hat T^{(12)}}(0,\frac{1}{2})$,  $\bold{x}_{\hat T^{(23)}}(\frac{1}{2},0)$, $\bold{x}_{\hat T^{(13)}}(\frac{1}{2},\frac{1}{2})$  and $c_{\hat T}  (\frac{1}{3},\frac{1}{3})$.
\begin{figure}[htbp!]
 \centering
 \caption{The reference triangle $ (\bold{x}_{\hat T^{(1)}},\bold{x}_{\hat T^{(2)}},\bold{x}_{\hat T^{(3)}})$.}
 \label{fig:referenceTriangle}
\end{figure}

Together with this assumption, we use lemma $3.2$ of \cite{LAMI} to obtain
\begin{equation}
\label{phtrinh24}
\int\limits_{\gamma _j^{(1)} } { N_{c_T }^b (\bold{x})\, \bold{n}_{\gamma _j^{(1)} } \text{d}\gamma (\bold{x})}  + \int\limits_{\gamma _k^{(1)} } { N_{c_T }^b (\bold{x})\, \bold{n}_{\gamma _k^{(1)} } \text{d}\gamma (\bold{x})} = J_T \left( {\hat \theta _{\hat \gamma _2^{(1)} } \hat{ \bold{n}_{\hat \gamma _2^{(1)} }}  + \hat \theta _{\hat \gamma _3^{(1)} } \hat  {\bold{n}_{\hat \gamma _3^{(1)} }} } \right)M_T^{ - 1}.
\end{equation}
\begin{equation}
\label{phtrinh25}
\int\limits_{V_{\bold{x}_i}  \cap T} {\nabla  N_{c_ T }^b (\bold{x})\text{d}\Omega}   = J_T \left( {\hat \theta _{\hat \gamma _2^{(2)} } \hat {\bold{n}_{\hat \gamma _2^{(2)} }}  + \hat \theta _{\hat \gamma _3^{(2)} } \hat  {\bold{n}_{\hat \gamma _3^{(2)} }} } \right)M_ T^{ - 1}.
\end{equation}
By directly computing the quantities on the reference element $\hat T$, we have
\begin{itemize}
\item The barycentric coordinates of a point $P(x^{(1)},x^{(2)})$ in the reference triangle $\hat T$ are $\hat \lambda_1(\bold{x}) = x^{(2)}$,  $\hat \lambda_2(\bold{x}) = 1 - x^{(1)} - x^{(2)}$  and  $\hat \lambda_3(\bold{x}) = x^{(1)}$ with $\bold{x}=(x^{(1)},x^{(2)})$. The basic cubic bubble function on the reference triangle $\hat T$ is  $N_{c_{\hat T}}^b (\bold{x}) = 27\lambda_1 (\bold{x})\lambda _2(\bold{x})\lambda_3(\bold{x})$.
\item The segments $\hat \gamma_1^{(1)}$, $\hat \gamma_1^{(2)}$ are  on the line $(d_1)~ x^{(2)} = -2x^{(1)} + 1$.
\item The segments $\hat \gamma_2^{(1)}$, $\hat \gamma_2^{(2)}$  are on the line $(d_2)~ x^{(2)} = x^{(1)}$.
\item The segments $\hat \gamma_3^{(1)}$, $\hat \gamma_3^{(2)}$  are on the line $(d_2)~ x^{(2)} = -0.5 x^{(1)} + 0.5$.
\item {The coefficients in (\ref{phtrinh24}) and (\ref{phtrinh25}) are computed by
\begin{equation}
\label{phtrinh26}
\hat\theta _{\hat \gamma _1^{(1)} }   =   \frac{{\sqrt 5 }}{{6}},\quad
\hat\theta _{\hat\gamma _1^{(2)} }   =\frac{{11\sqrt 5}}{{{\rm{96}}}},\quad
\hat\theta _{\hat \gamma _1^{(1)} }  = \frac{{16}}{{11}}\hat\theta _{\hat\gamma _1^{(2)} }.
\end{equation}
\begin{equation}
\label{phtrinh27}
\hat\theta _{\hat\gamma _2^{(1)} }   = \frac{{27\sqrt 2 }}{{162}},\quad
\hat \theta _{\hat \gamma _2^{(2)} }  =\frac{{ (11)(27)\sqrt 2}}{{{\rm{2592}}}}, \quad
\hat\theta _{\hat\gamma _2^{(1)} }  = \frac{{16}}{{11}}\hat \theta _{\hat \gamma _2^{(2)} }.
\end{equation}
\begin{equation}
\label{phtrinh28}
\hat \theta _{\hat \gamma _3^{(1)} } =  \frac{{16\sqrt {1.25} }}{{48}},\quad
\hat\theta _{\hat\gamma _3^{(2)} } = \frac{{11\sqrt {1.25} }}{{48}}, \quad
\hat \theta _{\hat \gamma _3^{(1)} }  = \frac{{16}}{{11}}\hat\theta _{\hat\gamma _3^{(2)} }.
\end{equation}}
\item The relationships between the normal vectors $\bold{n}_{\hat \gamma^{(1)}_i}$ and $\bold{n}_{\hat \gamma^{(2)}_i}$ with $i = \overline {1,3}$:
\begin{equation}
\label{phtrinh29}
\bold{n}_{\hat \gamma^{(1)}_1} = 2 \bold{n}_{\hat \gamma^{(2)}_1},\quad \bold{n}_{\hat \gamma^{(1)}_2} = 2 \bold{n}_{\hat \gamma^{(2)}_2}\quad \text{and} \quad \bold{n}_{\hat \gamma^{(1)}_3} = 2 \bold{n}_{\hat \gamma^{(2)}_3}.
\end{equation}
\end{itemize}
From (\ref{phtrinh24})-(\ref{phtrinh29}), we point out that
\begin{equation}
\label{phtrinh30}
\int\limits_{\Omega _{[\bold{x}_i ,\bold{x}_j ]}^s  \cap T} {\nabla  N_{c_T }^b (\bold{x})\text{d}\Omega}  + \int\limits_{\Omega _{[\bold{x}_k ,\bold{x}_i ]}^s  \cap T} {\nabla  N_{c_T }^b (\bold{x})\text{d}\Omega} = \frac{32}{11} \int\limits_{V_{\bold{x}_i}  \cap T} {\nabla N_{c_ T }^b (\bold{x})\text{d}\Omega}.
\end{equation}
Hence, we use the results of (\ref{phtrinh19}), (\ref{phtrinh23b}) and (\ref{phtrinh30}) to imply that
\begin{equation}
\label{phtrinh31}
\textit{ the coefficient of $p_{T^{(i)}} \bold{u}_{c_T}$ in $\overline b(\bold{u}_h,p_h)$}  = \frac{16}{11}. \textit{ the coefficient of $p_{T^{(i)}} \bold{u}_{c_T}$ in $b(\bold{u}_h,p_h)$}.
\end{equation}
With the computations of (\ref{phtrinh17}), (\ref{phtrinh18}) and (\ref{phtrinh31}), we conclude that
\begin{equation}
\label{phtrinh32}
 \int\limits_\Omega  {(\overline{\nabla \cdot\bold{b}_h} )q_h (\bold{x})\text{d}\Omega} = \frac{16}{11} \int\limits_\Omega  {(\nabla \cdot\bold{b}_h )q_h (\bold{x})\text{d}\Omega},
\end{equation}
Defining $\bold{u}^*_h = \boldsymbol{\ell}_h + \frac{11}{16}\bold{b}_h$, using (\ref{phtrinh32}) and the result of the first step, we get
\begin{equation}
\label{phtrinh33}
\int\limits_\Omega  {(\overline{\nabla \cdot\bold{u}^*_h} )q_h (\bold{x})\text{d}\Omega} = \int\limits_\Omega  {(\nabla \cdot\bold{u}_h )q_h (\bold{x})\text{d}\Omega}.
\end{equation}
Finally, due to the result of Theorem 3.1 in \cite{LAMI} and (\ref{phtrinh33}), the uniform inf-sup condition holds for the bilinear form $\overline b(\cdot,\cdot)$ on $\mathcal V_h^{\mathcal B} \times \mathcal V_h^{**}$.\\
\begin{itemize}
\renewcommand{\labelitemi}{$\blacksquare$}
\item \textit{ The hat bubble functions (\ref{hat})}\\
\end{itemize}
For each triangle $T \in \mathcal T_h$, the divergence of the hat bubble function is equal to a constant on each sub-triangle $\{T_{(i)}\}_{\overline{1,3}}$ of $T$, so we have
\begin{equation}
\label{phtrinh34}
\int\limits_{\Omega _{[\bold{x}_i ,\bold{x}_j ]}^s  \cap T} {\nabla  N_{c_T }^b (\bold{x})\text{d}\Omega}  + \int\limits_{\Omega _{[\bold{x}_k ,\bold{x}_i ]}^s  \cap T} {\nabla  N_{c_T }^b (\bold{x})\text{d}\Omega} = \frac{1}{2}\int\limits_{V_{\bold{x}_i}  \cap T} {\nabla N_{c_ T }^b (\bold{x})\text{d}\Omega}.
\end{equation}
By  (\ref{phtrinh19}), (\ref{phtrinh23b}) and  (\ref{phtrinh34}),  we obtain
\begin{equation}
\label{phtrinh35}
\textit{ the coefficient of $p_{T^{(i)}} \bold{u}_{c_T}$ in $\overline b(\bold{u}_h,p_h)$}  =  \textit{ the coefficient of $p_{T^{(i)}} \bold{u}_{c_T}$ in $b(\bold{u}_h,p_h)$},
\end{equation}
which implies that
\begin{equation}
\label{phtrinh36}
\int\limits_\Omega  {(\overline{\nabla \cdot\bold{b}_h} )q_h (\bold{x})\text{d}\Omega} = \int\limits_\Omega  {(\nabla \cdot\bold{b}_h )q_h (\bold{x})\text{d}\Omega}.
\end{equation}
Therefore,
\begin{equation}
\label{phtrinh37}
\int\limits_\Omega  {(\overline{\nabla \cdot\bold{u}_h} )q_h (\bold{x})\text{d}\Omega} = \int\limits_\Omega  {(\nabla \cdot\bold{u}_h )q_h (\bold{x})\text{d}\Omega}.
\end{equation}
\newline
In three dimensions,  we also compute the coefficient of  $q_{{T^{(i)} }} \bold{u}_{c_T}$ (\ref{phtrinh23}) on the following tetrahedron $T$ constructed from four vertices  $\{\bold{x}_{T^{(i)}},\bold{x}_{T^{(j)}},\bold{x}_{T^{(k)}},\bold{x}_{T^{(l)}}\}$,
\begin{figure}[htbp!]
  \centering
  \caption{A tetrahedron $(\bold{x}_{T^{(i)}},\bold{x}_{T^{(j)}},\bold{x}_{T^{(k)}},\bold{x}_{T^{(l)}})$ belonging to $\mathcal T_h$.}
  \label{fig:tetrahedron}
  \end{figure}
where $c_T$ is the centroid  of $T$, $\bold{x}_f^{(ijk)}$ is the barycentric point of the triangular face $(\bold{x}_{T^{(i)}},\bold{x}_{T^{(j)}},\bold{x}_{T^{(k)}})$, $\bold{x}_e^{(ij)}$ is the midpoint of the edge $[\bold{x}_{T^{(i)}},\bold{x}_{T^{(j)}}]$,  $i, j, k, l$ belong to $\{1,2,3,4\}$.\\
In this particular case, the coefficient of  $q_{{T^{(i)} }} \bold{u}_{c_T}$ (\ref{phtrinh23}) of $\overline b(\bold{u},q)$ is computed as
\begin{eqnarray}
\label{phtrinh38}
\frac{1}{3}\left( {\int\limits_{\Omega _{\left[ {\bold{x}_{T^{(i)} } ,\bold{x}_{T^{(j)} } } \right]}^s \cap T} {\nabla N_{c_T }^b (\bold{x})\text{d}\Omega}  + \int\limits_{\Omega _{\left[ {\bold{x}_{T^{(i)} } ,\bold{x}_{T^{(k)} } } \right]}^s  \cap T} {\nabla N_{c_T }^b (\bold{x})\text{d}\Omega}  + \int\limits_{\Omega _{\left[ {\bold{x}_{T^{(i)} } ,\bold{x}_{T^{(l)} } } \right]}^s  \cap T} {\nabla N_{c_T }^b (\bold{x})\text{d}\Omega} } \right) &&\nonumber\\
\nonumber\\
=\frac{1}{3}\left( \begin{array}{l}
 \int\limits_{\left( {c_T ,\bold{x}_f^{(ijk)} ,\bold{x}_{T^{(j)} } } \right)} {N_{c_T }^b (\bold{x})\bold{n}_{\left( {c_T ,\bold{x}_f^{(ijk)} ,\bold{x}_{T^{(j)} } } \right)} \text{d}\gamma (\bold{x})}  + \int\limits_{\left( {c_T ,\bold{x}_f^{(ijl)} ,\bold{x}_{T^{(j)} } } \right)} {N_{c_T }^b (\bold{x})\bold{n}_{\left( {c_T ,\bold{x}_f^{(ijl)} ,\bold{x}_{T^{(j)} } } \right)} \text{d}\gamma (\bold{x})} +  \\
\\
 \int\limits_{\left( {c_T ,\bold{x}_f^{(ijk)} ,\bold{x}_{T^{(k)} } } \right)} {N_{c_T }^b (\bold{x})\bold{n}_{\left( {c_T ,\bold{x}_f^{(ijk)} ,\bold{x}_{T^{(k)} } } \right)} \text{d}\gamma (\bold{x})}  + \int\limits_{\left( {c_T ,\bold{x}_f^{(ikl)} ,\bold{x}_{T^{(k)} } } \right)} {N_{c_T }^b (\bold{x})\bold{n}_{\left( {c_T ,\bold{x}_f^{(ikl)} ,\bold{x}_{T^{(k)} } } \right)} \text{d}\gamma (\bold{x})} + \\
\\
 \int\limits_{\left( {c_T ,\bold{x}_f^{(ijl)} ,\bold{x}_{T^{(l)} } } \right)} {N_{c_T }^b (\bold{x})\bold{n}_{\left( {c_T ,\bold{x}_f^{(ijl)} ,\bold{x}_{T^{(l)} } } \right)} \text{d}\gamma (\bold{x})}  + \int\limits_{\left( {c_T ,\bold{x}_f^{(ikl)} ,\bold{x}_{T^{(l)} } } \right)} {N_{c_T }^b (\bold{x})\bold{n}_{\left( {c_T ,\bold{x}_f^{(ikl)} ,\bold{x}_{T^{(l)} } } \right)} \text{d}\gamma (\bold{x})}\\
 \end{array} \right),
\end{eqnarray}
where vectors $\bold{n}_{\left( {c_T ,\bold{x}_f^{(ijk)} ,\bold{x}_{T^{(j)} } } \right)}$, $\bold{n}_{\left( {c_T ,\bold{x}_f^{(ijl)} ,\bold{x}_{T^{(j)} } } \right)}$, $\bold{n}_{\left( {c_T ,\bold{x}_f^{(ijk)} ,\bold{x}_{T^{(k)} } } \right)}$, $\bold{n}_{\left( {c_T ,\bold{x}_f^{(ikl)} ,\bold{x}_{T^{(k)} } } \right)}$, $\bold{n}_{\left( {c_T ,\bold{x}_f^{(ijl)} ,\bold{x}_{T^{(l)} } } \right)}$, $\bold{n}_{\left( {c_T ,\bold{x}_f^{(ikl)} ,\bold{x}_{T^{(l)} } } \right)}$ whose length is equal to  measure of triangular faces ${\left( {c_T ,\bold{x}_f^{(ijk)} ,\bold{x}_{T^{(j)} } } \right)}$, \\${\left( {c_T ,\bold{x}_f^{(ijl)} ,\bold{x}_{T^{(j)} } } \right)}$, ${\left( {c_T ,\bold{x}_f^{(ijk)} ,\bold{x}_{T^{(k)} } } \right)}$, ${\left( {c_T ,\bold{x}_f^{(ikl)} ,\bold{x}_{T^{(k)} } } \right)}$,  ${\left( {c_T ,\bold{x}_f^{(ijl)} ,\bold{x}_{T^{(l)} } } \right)}$ and ${\left( {c_T ,\bold{x}_f^{(ikl)} ,\bold{x}_{T^{(l)} } } \right)}$, are the outward normal vectors of $T\cap \Omega^s_{[\bold{x}_i,\bold{x}_j]}$, $T\
cap \Omega^s_{[\bold{x}_i,\bold{x}_k]}$ and $T\cap \Omega^s_{[\bold{x}_i,\bold{x}_l]}$.\\
\newline
We also get the coefficient of  $q_{{T^{(i)} }} \bold{u}_{c_T}$ of $ b(u,q)$, as follows:
\begin{equation}
\label{phtrinh39}
\left( \begin{array}{l}
 \int\limits_{\left( {c_T ,\bold{x}_f^{(ijk)} ,\bold{x}_e^{(i,j)} } \right)} {N_{c_T }^b (\bold{x})\bold{n}_{\left( {c_T ,\bold{x}_f^{(ijk)} ,\bold{x}_e^{(ij)} } \right)} \text{d}\gamma (\bold{x})}  + \int\limits_{\left( {c_T ,\bold{x}_f^{(ijl)} ,\bold{x}_e^{(ij)} } \right)} {N_{c_T }^b (\bold{x})\bold{n}_{\left( {c_T ,\bold{x}_f^{(ijl)} ,\bold{x}_e^{(ij)} } \right)} \text{d}\gamma (\bold{x})} + \\
\\
 \int\limits_{\left( {c_T ,\bold{x}_f^{(ijk)} ,\bold{x}_e^{(i k)} } \right)} {N_{c_T }^b (\bold{x})\bold{n}_{\left( {c_T ,\bold{x}_f^{(ijk)} ,\bold{x}_e^{(i k)} } \right)} \text{d}\gamma (\bold{x})}  + \int\limits_{\left( {c_T ,\bold{x}_f^{(ikl)} ,\bold{x}_e^{(i k)} } \right)} {N_{c_T }^b (\bold{x})\bold{n}_{\left( {c_T ,\bold{x}_f^{(ikl)} ,\bold{x}_e^{(i k)} } \right)} \text{d}\gamma (\bold{x})}  +\\
\\
 \int\limits_{\left( {c_T ,\bold{x}_f^{(ijl)} ,\bold{x}_e^{(i,l)} } \right)} {N_{c_T }^b (\bold{x})\bold{n}_{\left( {c_T ,\bold{x}_f^{(ijl)} ,\bold{x}_e^{(i l)} } \right)} \text{d}\gamma (\bold{x})}  + \int\limits_{\left( {c_T ,\bold{x}_f^{(ikl)} ,\bold{x}_e^{(il)} } \right)} {N_{c_T }^b (\bold{x})\bold{n}_{\left( {c_T ,\bold{x}_f^{(ikl)} ,\bold{x}_e^{(i l)} } \right)} \text{d}\gamma (\bold{x})} \\
 \end{array} \right).
\end{equation}
Furthermore, we have relationships between normal vectors in the two formulas (\ref{phtrinh38}) and (\ref{phtrinh39})
\begin{eqnarray}
\label{phtrinh40}
\bold{n}_{\left( {c_T ,\bold{x}_f^{(ijk)} ,\bold{x}_{T^{(j)} } } \right)}  &=& 2\bold{n}_{\left( {c_T ,\bold{x}_f^{(ijk)} ,\bold{x}_e^{(i j)} } \right)} ,\,\bold{n}_{\left( {c_T ,\bold{x}_f^{(ijl)} ,\bold{x}_{T^{(j)} } } \right)}  = 2\bold{n}_{\left( {c_T ,\bold{x}_f^{(ijl)} ,\bold{x}_e^{(i j)} } \right)}, \nonumber\\
\nonumber \\
\bold{n}_{\left( {c_T ,\bold{x}_f^{(ijk)} ,\bold{x}_{T^{(k)} } } \right)}  &=& 2\bold{n}_{\left( {c_T ,\bold{x}_f^{(ijk)} ,\bold{x}_e^{(i k)} } \right)}, \bold{n}_{\left( {c_T ,\bold{x}_f^{(ikl)} ,\bold{x}_{T^{(k)} } } \right)}  = 2\bold{n}_{\left( {c_T ,\bold{x}_f^{(ikl)} ,\bold{x}_e^{(i k)} } \right)}, \nonumber \\
\nonumber \\
\bold{n}_{\left( {c_T ,\bold{x}_f^{(ijl)} ,\bold{x}_{T^{(l)} } } \right)}  &=& 2\bold{n}_{\left( {c_T ,\bold{x}_f^{(ijl)} ,\bold{x}_e^{(il)} } \right)} ,\bold{n}_{\left( {c_T ,\bold{x}_f^{(ikl)} ,\bold{x}_{T^{(l)} } } \right)}  = 2\bold{n}_{\left( {c_T ,\bold{x}_f^{(ikl)} ,\bold{x}_e^{(i l)} } \right)}.
\end{eqnarray}
\newpage
\noindent{\bf For the bFS-FEM method}\\
\newline
In a similar manner to the calculations for bES-FEM, the coefficient of $q_{T^{(i)}}\bold{u}_{c_T}$ is equal to
\begin{eqnarray}
&&\sum\limits_{\scriptstyle \,f_{T^{(i)} }  \in \mathcal F_{T^{(i)} } ,\Omega _{f_{T^{(i)} } }^s  \in \mathcal T_h^{**}  \hfill \atop
  \scriptstyle \,\,K \in \mathcal T_h ,~ \mathcal F_K  \cap \mathcal F_T  = f_{T^{(i)} }  \hfill} {\left[ \begin{array}{l}
 \frac{{m\left( {V_{\bold{x}_{T^{(i)} } }  \cap T \cap \Omega _{f_{T^{(i)} } }^s } \right)}}{{m\left( {\Omega _{f_{T^{(i)} } }^s } \right)}} \\
 \frac{{m\left( {V_{\bold{x}_{T^{(i)} } }  \cap K \cap \Omega _{f_{T^{(i)} } }^s } \right)}}{{m\left( {\Omega _{f_{T^{(i)} } }^s } \right)}} \\
 \end{array} \right]} \,\,\int\limits_{\Omega _{f_{T^{(i)} } }^s  \cap T} {\nabla N_{c_T }^b (\bold{x})\text{d}\Omega} \nonumber \\
&=&\frac{1}{d}\sum\limits_{\,f_{T^{(i)} }  \in \mathcal F_{T^{(i)} } ,\Omega _{f_{T^{(i)} } }^s  \in \mathcal T_h^{**} }~~{\int\limits_{\Omega _{f_{T^{(i)} } }^s  \cap T} {\nabla N_{c_T }^b (\bold{x})\text{d}\Omega} }.
\end{eqnarray}
With a tetrahedral $T = (\bold{x}_{T^{(i)}},\bold{x}_{T^{(j)}},\bold{x}_{T^{(k)}},\bold{x}_{T^{(l)}}) \in \mathcal T_h$ (see Figure~\ref{fig:tetrahedron}),
we obtain the following coefficient of $q_{T^{(i)}}\bold{u}_{c_T}$ of $\overline b(\bold{u},q)$ for bFS-FEM
\begin{eqnarray*}
\label{bFS1}
\int\limits_{\Omega _{\left(\bold{x}_{T^{(i)}}  ,\bold{x}_{T^{(j)}},\bold{x}_{T^{(k)}}\right)}^s  \cap T} {\nabla N_{c_T }^b (\bold{x})\text{d}\Omega}  + \int\limits_{\Omega _{\left(\bold{x}_{T^{(i)}} ,\bold{x}_{T^{(j)}}  ,\bold{x}_{T^{(l)}} \right)}^s  \cap T} {\nabla N_{c_T }^b (\bold{x})\text{d}\Omega}  + \int\limits_{\Omega _{(\bold{x}_{T^{(i)}} ,\bold{x}_{T^{(k)}} ,\bold{x}_{T^{(l)}} )}^s  \cap T} {\nabla N_{c_T }^b (\bold{x})\text{d}\Omega} 
\end{eqnarray*}
\begin{eqnarray}
= \int\limits_{(c_T ,\bold{x}_{T^{(j)} } ,\bold{x}_{T^{(k)} } )} {N_{c_T }^b (\bold{x})\bold{n}_{(c_T ,\bold{x}_{T^{(j)} } ,\bold{x}_{T^{(k)} } )} \text{d}\gamma (\bold{x})}  + \int\limits_{(c_T ,\bold{x}_{T^{(j)} } ,\bold{x}_{T^{(l)} } )} {N_{c_T }^b (\bold{x})\bold{n}_{(c_T ,\bold{x}_{T^{(j)} } ,\bold{x}_{T^{(l)} } )} \text{d}\gamma (\bold{x})}  +&&\nonumber \\
\int\limits_{(c_T ,\bold{x}_{T^{(k)} } ,\bold{x}_{T^{(l)} } )} {N_{c_T }^b (\bold{x})\bold{n}_{(c_T ,\bold{x}_{T^{(k)} } ,\bold{x}_{T^{(l)} } )} \text{d}\gamma (\bold{x})},
\end{eqnarray}
where normal unit vectors $\bold{n}_{(c_T ,\bold{x}_{T^{(j)} } ,\bold{x}_{T^{(k)} } )}$, $\bold{n}_{(c_T ,\bold{x}_{T^{(j)} } ,\bold{x}_{T^{(l)} } )}$ and $\bold{n}_{(c_T ,\bold{x}_{T^{(k)} } ,\bold{x}_{T^{(l)} } )}$ of the tetrahedron ${(c_T ,\bold{x}_{T^{(i)}},\bold{x}_{T^{(j)} } ,\bold{x}_{T^{(k)} } )}$, ${(c_T ,\bold{x}_{T^{(i)}},\bold{x}_{T^{(j)} } ,\bold{x}_{T^{(l)} } )}$ and ${(c_T ,\bold{x}_{T^{(i)}},\bold{x}_{T^{(k)} } ,\bold{x}_{T^{(l)} } )}$ are measured by the area of triangular faces ${(c_T ,\bold{x}_{T^{(j)} } ,\bold{x}_{T^{(k)} } )}$, ${(c_T ,\bold{x}_{T^{(j)} } ,\bold{x}_{T^{(l)} } )}$ and ${(c_T ,\bold{x}_{T^{(k)} } ,\bold{x}_{T^{(l)} } )}$, respectively. Additionally, normal vectors in (\ref{phtrinh38}) and (\ref{bFS1}) relate  together
\begin{eqnarray}
\label{bFS2b}
&&\bold{n}_{\left( {c_T ,\bold{x}_{T^{(j)} } ,\bold{x}_f^{(ijl)} } \right)}  = \bold{n}_{\left( {c_T ,\bold{x}_{T^{(k)} } ,\bold{x}_f^{(ikl)} } \right)}  = \frac{1}{3}\bold{n}_{\left( {c_T ,\bold{x}_{T^{(k)} } ,\bold{x}_{T^{(j)} } } \right)} \nonumber\\
&&\bold{n}_{\left( {c_T ,\bold{x}_{T^{(j)} } ,\bold{x}_f^{(ijk)} } \right)}  = \bold{n}_{\left( {c_T ,\bold{x}_{T^{(l)} } ,\bold{x}_f^{(ikl)} } \right)}  = \frac{1}{3}\bold{n}_{\left( {c_T ,\bold{x}_{T^{(j)} } ,\bold{x}_{T^{(l)} } } \right)}, \nonumber \\
&&\bold{n}_{\left( {c_T ,\bold{x}_{T^{(k)} } ,\bold{x}_f^{(ijk)} } \right)}  = \bold{n}_{\left( {c_T ,\bold{x}_{T^{(l)} } ,\bold{x}_f^{(ijl)} } \right)}  = \frac{1}{3}\bold{n}_{\left( {c_T ,\bold{x}_{T^{(k)} } ,\bold{x}_{T^{(l)} } } \right)}.
\end{eqnarray}
From  (\ref{phtrinh38})-(\ref{bFS2b}), there exist the two positive constants  $\alpha_1$, $\alpha_2$ satisfying
\begin{itemize}
\item \textit{ the coefficient of $p_{T^{(i)}} \bold{u}_{c_T}$ in $\overline b(\bold{u}_h,p_h)$ of the 3D bES-FEM = $\alpha_1$. the coefficient of $p_{T^{(i)}} \bold{u}_{c_T}$ in $b(\bold{u}_h,p_h)$ and}
\item \textit{ the coefficient of $p_{T^{(i)}} \bold{u}_{c_T}$ in $\overline b(\bold{u}_h,p_h)$ of the bFS-FEM = $\alpha_2$. the coefficient of $p_{T^{(i)}} \bold{u}_{c_T}$ in $\overline b(\bold{u}_h,p_h)$ of the 3D bES-FEM},
\end{itemize}
which lead to
\begin{eqnarray}
\label{phtrinh42}
&&\overline b(\bold{u}_h,q_h) ~\text{of the 3D bES-FEM} = \alpha_1 \int\limits_\Omega  {(\nabla \cdot\bold{b}_h )q_h (\bold{x})\text{d}\Omega},\nonumber\\
&&\overline b(\bold{u}_h,q_h) ~\text{of the bFS-FEM} = \alpha_2 ~\overline b(\bold{u}_h,q_h)~\text{of the 3D bES-FEM  implying that}\nonumber\\
&&\overline b(\bold{u}_h,q_h) ~\text{of the bFS-FEM} = \alpha_1 \alpha_2 \int\limits_\Omega  {(\nabla \cdot\bold{b}_h )q_h (\bold{x})\text{d}\Omega}.
\end{eqnarray}
Therefore, for each the 2D/3D bES-FEM or the bFS-FEM, we choose the coefficient $\alpha$ that is equal to $\alpha_1$ or $\alpha_1\alpha_2$, respectively. $\hspace{9.32cm}\square$\\
\newline
From the results of the two lemmas \ref{mathproper}.1 and \ref{mathproper}.2, we deduce that there are two positive constants $\alpha_3$, $\alpha_4$ depending on $\alpha_1$, $\alpha_2$, such that
\begin{eqnarray}
\label{phtrinh44}
\text{the 2D/3D bES-FEM method} \quad \quad \int\limits_\Omega  {(\overline{\nabla \cdot\bold{u}^*_h} )q_h (\bold{x})\text{d}\Omega} &=& \int\limits_\Omega  {(\nabla \cdot\bold{u}_h )q_h (\bold{x})\text{d}\Omega}. \nonumber \\
\text{the bFS-FEM method} \quad \quad \int\limits_\Omega  {(\overline{\nabla \cdot\bold{u}^{**}_h} )q_h (\bold{x})\text{d}\Omega} &=& \int\limits_\Omega  {(\nabla \cdot\bold{u}_h )q_h (\bold{x})\text{d}\Omega}.
\end{eqnarray}
with $\bold{u}_h = \boldsymbol{\ell}_h + \bold{b}_h$, $\bold{u}^*_h =  \boldsymbol{\ell}_h + \alpha_3 \bold{b}_h$ and $\bold{u}^{**}_h =  \boldsymbol{\ell}_h + \alpha_4 \bold{b}_h$ in $\Omega$.\\
Hence theorem \ref{mathproper}.2 is proven.
$\hspace{10.5cm}\square$\\
\newline
Additionally, the bilinear form $c(\cdot,\cdot)$ is continuous, symmetric and positive semi-definite, {\it i.e.}
\begin{displaymath}
\overline c(q,q) \ge 0, \quad q \in L^2_0(\Omega).
\end{displaymath}
\newline
{\bf Theorem \ref {mathproper}.3} \textit{(Convergence)}\\
\textit{We assume that $(\bold u,p)$ and $(\bold u_h,p_h)$ are the two pair solutions of the problems (\ref{phtrinh05},\ref{phtrinh05.1}) and (\ref{phtrinh4}a,\ref{phtrinh4}b), then we get the following error estimation}
\begin{equation}
\label{conver1}
||\bold u- \bold u_h||_{\mathcal V^{\mathcal B}_h} + ||p-p_h||_{L^2(\Omega)} \le C \left(\mathop {\inf}\limits_{{{\bf{w}}_h} \in \mathcal V_h^{\mathcal B}} \{||\bold w_h - \bold u||_{(\mathbb H^1(\Omega))^2} \}  +  \mathop {\inf }\limits_{{q_h} \in \mathcal V_h^{**}} \{||q_h - p||_{L^2(\Omega)}\} \right )+ \mathcal O (h)
\end{equation}
where $C$ is a positive constant and independent on $h$. This coefficient $h$ is defined by
\begin{equation}
\label{conver1a}
h =  \max \left\{ \mathop {\sup}\limits_{ K^* \in \mathcal M^*} \text{diam} (K^*), \mathop {\sup}\limits_{ K^{**} \in \mathcal M^{**}} \text{diam} (K^{**})  \right\},
\end{equation}
and a radius of the circumscribed circle for each element $K^*$ of $\mathcal M^*$, $K^{**}$ of $\mathcal M^{**}$ is denoted by ``diam($K^*$)", ``diam($K^{**}$)", respectively.\\
\newline
\textit{Proof:}\quad Let us consider any $\bold w_h \in \mathcal V^{\mathcal B}_h(\lambda)$ defined by
\begin{displaymath}
\mathcal V^{\mathcal B}_h(\lambda) = \left\{ \bold w_h \in \mathcal V^{\mathcal B}_h ~|~  \overline b(\bold w_h, q_h) = \frac{1}{\lambda} \overline c(p_h,q_h), \forall q_h\in \mathcal V^{**}_h \right \}
\end{displaymath}
This implies $\overline b(\bold w_h - \bold u_h, q_h) = 0$, for all $q_h\in \mathcal V^{**}_h$, \textit{i.e}, $\bold w_h - \bold u_h$  is an element of $\mathcal V^{\mathcal B}_{h,0}\subset [H^1(\Omega)]^d$. Then, by applying the coercivity (\ref{phtrinh5a}), one has
\begin{eqnarray}
\label{conver2}
\alpha_0 ||\bold w_h - \bold u_h||^2 &\le&  \overline a(\bold w_h - \bold u_h,\bold w_h - \bold u_h) \nonumber \\
&=& \left[\overline a(\bold w_h,\bold w_h - \bold u_h) - a(\bold u,\bold w_h - \bold u_h) +  a(\bold u,\bold w_h - \bold u_h) - \overline a(\bold u_h,\bold w_h - \bold u_h) \right] \nonumber \\
&=&  \left[ \begin{array}{l} \overline a(\bold w_h - \bold u,\bold w_h - \bold u_h) + \overline a(\bold u, \bold w_h - \bold u_h) - a(\bold u, \bold w_h - \bold u_h) -\\
b(\bold w_h - \bold u_h,p - q_h) + b(\bold u_h- \bold w_h,q_h) - \overline b(\bold u_h - \bold w_h,p^{\mathcal M^{**}}_h) \end{array} \right],
\end{eqnarray}
where $p^{\mathcal M^{**}_h} \in \mathcal V_h^{**}$ is a characteristic function defined by
\begin{displaymath}
p^{\mathcal M^{**}}_h|_{K^{**}} = \frac{1}{\text{m}(K^{**})} \int\limits_{K^{**}} {p(\bold x) \text{d} \Omega},
\end{displaymath}
Note that we have $a(\bold u,\bold w_h - \bold u_h) - \overline a(\bold u_h,\bold w_h - \bold u_h) = b(\bold w_h - \bold u_h,p)$, this is a result of
$(\ref{phtrinh4}a)$ subtracted from $(\ref{phtrinh05})$.\\
Inequality (\ref{conver2}) continues to be evaluated as follows
\begin{eqnarray}
\label{conver3}
\alpha_0||\bold w_h - \bold u_h ||^2_{\mathcal V^{\mathcal B}_h} &\le&
\left[ \begin{array}{l}
\vspace{0.5cm}
\underbrace {\alpha_1 ||\bold w_h - \bold u_h||_{\mathcal V^{\mathcal B}_h} ||\bold w_h - \bold u||_{\mathcal V^{\mathcal B}_h}}_{\ge \overline a(\bold w_h - \bold u,\bold w_h - \bold u_h)} + \underbrace {||\nabla.(\bold w_h -\bold u_h)||_{(L^2(\Omega))^2}~||q_h - p||_{L^2(\Omega)}}_{\ge b(\bold w_h - \bold u_h,p - q_h)} + \\
\vspace{0.4cm}
 \underbrace{\lambda^{\bold D}_{\max}~ ||\boldsymbol{\varepsilon}(\bold u) - \overline {\boldsymbol{\varepsilon} (\bold u)}||_{(L^2(\Omega))^2} ||\boldsymbol{\varepsilon}(\bold w_h - \bold u_h)||_{(L^2(\Omega))^2}}_{\ge \overline a(\bold u,\bold w_h - \bold u_h) - a(\bold u,\bold w_h - \bold u_h)} + \\
 \vspace{0.4cm}
b(\bold u_h- \bold w_h,q_h) - b(\bold u_h- \bold w_h,p^{\mathcal M^*}_h)+b(\bold u_h- \bold w_h,p^{\mathcal M^*}_h)- \overline b(\bold u_h- \bold w_h,p)+\\
\overline b(\bold u_h- \bold w_h,p)- \overline b(\bold u_h - \bold w_h,p^{\mathcal M^{**}}_h)
\end{array} \right], \nonumber\\
\end{eqnarray}
where eigenvalues of the material matrix $\bold{D}$  are upper bounded by $\lambda^{\bold D}_{\max}$, and $p^{\mathcal M^*}_h$ is a characteristic function defined by
\begin{displaymath}
p^{\mathcal M^*}_h|_{K^*} = \frac{1}{\text{m}(K^*)} \int\limits_{K^*} {p(\bold x) \text{d} \Omega},
\end{displaymath}
with the pressure solution $p$ of (\ref{phtrinh05}, \ref{phtrinh05.1}), on each element $K^* \in \mathcal M^{*}$.\\
Besides, we have the following estimations
\begin{equation*}
\label{conver3a}
\frac{b(\bold u_h - \bold w_h,q_h) - b(\bold u_h - \bold w_h,p^{\mathcal M^*}_h)}{||\bold w_h - \bold u_h||_{\mathcal V^{\mathcal B}_h}} \le \frac{||\bold w_h - \bold u_h||_{(L^2(\Omega))^2}}{||\bold w_h - \bold u_h||_{\mathcal V^{\mathcal B}_h}}||q_h - p||_{L^2(\Omega)} + \frac{||\bold w_h - \bold u_h||_{(L^2(\Omega))^2}}{||\bold w_h - \bold u_h||_{\mathcal V^{\mathcal B}_h}}||p^{\mathcal M^*}_h - p||_{L^2(\Omega)},
\end{equation*}
\begin{equation}
\label{conver3aa}
\frac{\overline b(\bold u_h - \bold w_h,p) - \overline b(\bold u_h - \bold w_h,p^{\mathcal M^{**}}_h)}{||\bold w_h - \bold u_h||_{\mathcal V^{\mathcal B}_h}}  \mathop {\le}\limits_{\text{Holder}} \frac{||\overline{\nabla.(\bold w_h - \bold u_h)}||_{L^2(\Omega)}}{||\bold w_h - \bold u_h||_{\mathcal V^{\mathcal B}_h}}||p^{\mathcal M^{**}}_h - p||_{L^2(\Omega)}
\end{equation}
because of (26), (27, (30), (31) and (45) found in \cite{GRL}. Moreover, we have
\begin{eqnarray}
\label{conver3b}
&&\frac{b(\bold u_h- \bold w_h,p^{\mathcal M^*}_h)- \overline b(\bold u_h - \bold w_h,p)}{||\bold w_h - \bold u_h||_{\mathcal V^{\mathcal B}_h}}  \nonumber\\
&=&\sum\limits_{{K^*} \in {\mathcal M^*}} {\left[ \begin{array}{l}
\vspace{0.4cm}
\frac{1}{{\text{m}({K^*})}}\left( {\int\limits_{{K^*}} {p(\bold x)} {\rm{d}}\Omega } \right)\left( {\int\limits_{{K^*}} {\frac{{\nabla .{(\bold u_h- \bold w_h)}(\bold x)}}{{||{\bold u_h- \bold w_h}|{|_{\mathcal V_h^{\mathcal B}}}}}} {\text{d}}\Omega } \right) - \\
\frac{1}{{\text{m}({K^*})}}\left( {\int\limits_{{K^*}} {\frac{{\nabla .{(\bold u_h- \bold w_h)}(\bold x)}}{{||{\bold u_h- \bold w_h}|{|_{\mathcal V_h^{\mathcal B}}}}}} {\rm{d}}\Omega } \right)\left( {\int\limits_{{K^*}} {p(\bold x)} {\text{d}}\Omega } \right)
\end{array} \right]} = 0.
\end{eqnarray}
Using (\ref{conver3a})-(\ref{conver3b}), the inequality (\ref{conver3}) is rewritten as follows
\begin{eqnarray}
\label{conver3c}
&&\left[ \begin{array}{l}
\vspace{0.5cm}
\alpha_1 ||\bold w_h - \bold u||_{\mathcal V^{\mathcal B}_h}  + \frac{||\nabla.(\bold w_h -\bold u_h)||_{L^2(\Omega)}}{||\bold w_h - \bold u_h||_{\mathcal V^{\mathcal B}_h}}~||q_h - p||_{L^2(\Omega)} + \\
\vspace{0.5cm}
 \lambda^{\bold D}_{\max}~ ||\boldsymbol{\varepsilon}(\bold u) - \overline {\boldsymbol{\varepsilon} (\bold u)}||_{(L^2(\Omega))^2} \frac{||\boldsymbol{\varepsilon}(\bold w_h - \bold u_h)||_{(L^2(\Omega))^2}}{||\bold w_h - \bold u_h||_{\mathcal V^{\mathcal B}_h}} \\
 \vspace{0.5cm}
+\frac{||\bold w_h - \bold u_h||_{(L^2(\Omega))^2}}{||\bold w_h - \bold u_h||_{\mathcal V^{\mathcal B}_h}}||q_h - p||_{L^2(\Omega)} + \frac{||\bold w_h - \bold u_h||_{(L^2(\Omega))^2}}{||\bold w_h - \bold u_h||_{\mathcal V^{\mathcal B}_h}}||p^{\mathcal M^*}_h - p||_{L^2(\Omega)} \\
\frac{||\overline{\nabla.(\bold w_h - \bold u_h)}||_{L^2(\Omega)}}{||\bold w_h - \bold u_h||_{\mathcal V^{\mathcal B}_h}}||p^{\mathcal M^{**}}_h - p||_{L^2(\Omega)}
\end{array} \right] \nonumber \\
&\ge& \alpha_0||\bold w_h - \bold u_h ||_{\mathcal V^{\mathcal B}_h} \ge \alpha_0(||\bold u_h - \bold u||_{\mathcal V^{\mathcal B}_h} - ||\bold w_h - \bold u ||_{\mathcal V^{\mathcal B}_h})
\end{eqnarray}
Let us subtract (\ref{phtrinh05}) from (\ref{phtrinh4}a), getting
\begin{equation}
\label{conver4}
\overline b(\bold v_h, p_h) - b(\bold v_h, p) =  a(\bold u, \bold v_h) - \overline a(\bold u_h, \bold v_h), \quad \forall \bold v_h \in \mathcal V^{\mathcal B}_h,
\end{equation}
so that for $q_h \in \mathcal V^{**}_h$, it follows
\begin{eqnarray}
\label{conver5}
\overline b(\bold v_h, p_h - q_h) 
&=& a(\bold u, \bold v_h) - \overline a(\bold u_h, \bold v_h) +  b(\bold v_h, p) - \overline b(\bold v_h, q_h)
\end{eqnarray}
Transforming $\overline b(\bold v_h, p_h - q_h)$ in the stability property (\ref{phtrinh5}) by (\ref{conver5}), we have
\begin{equation}
\label{conver6}
\mathop {\sup}\limits_{\bold v_h \in \mathcal V^{\mathcal B}_h, \bold v_h \ne 0} \frac{a(\bold u, \bold v_h) - \overline a(\bold u_h, \bold v_h) +  b(\bold v_h, p) - \overline b(\bold v_h, q_h)}{||\bold v_h||_{\mathcal V^{\mathcal B}_h}}=\mathop {\sup}\limits_{\bold v_h \in \mathcal V^{\mathcal B}_h, \bold v_h \ne 0} \frac{\overline b(\bold v_h,p_h-q_h)}{||\bold v_h||_{\mathcal V^{\mathcal B}_h}} \ge \beta_0 ||p_h - q_h||_{L^2(\Omega)}.
\end{equation}
Now, we estimate each part in the left hand side of (\ref{conver6}):
\begin{eqnarray}
\label{conver7}
&&\frac{b(\bold v_h,p) - \overline b(\bold v_h,q_h)}{||\bold v||_{\mathcal V^{\mathcal B}_h}} = \frac{b(\bold v_h, p)-b(\bold v_h, p^{\mathcal M^*}_h) + b(\bold v_h, p^{\mathcal M^*}_h)  - \overline b(\bold v_h,p) + \overline b(\bold v_h,p) - \overline b(\bold v_h, q_h)}{{||\bold v_h||}_{\mathcal V^{\mathcal B}_h}} \nonumber \\
&\mathop {\le}\limits_{\text{Holder}} & \frac{||\nabla.\bold v_h||_{L^2(\Omega)}}{{{||{{\bf{v}}_h}|{|_{\mathcal V_h^{\mathcal B}}}}}} ||p-p^{\mathcal M^*}_h||_{L^2(\Omega)} + \frac{||\overline{\nabla. \bold v_h}||_{L^2(\Omega)}}{||\bold v_h||_{\mathcal V^{\mathcal B}_h}}||p - q_h||_{L^2(\Omega)}
\end{eqnarray}
because of
 $\int\limits_\Omega  {\frac{{ (\nabla.\bold v_h) p^{\mathcal M^*}_h - (\overline {\nabla .{{\bold{v}}_h}}) p}}{{||{{\bold{v}}_h}|{|_{\mathcal V_h^{\mathcal B}}}}}}\,{\rm{d}}\Omega = 0$, explained as (\ref{conver3b}).\\
\newline
For the other part of (\ref{conver6}), thanks to two equations (63a) and (63b) of \cite{LiHu}, one writes
\begin{eqnarray}
\label{conver8}
&&\frac{a(\bold u, \bold v_h) - \overline{a}(\bold u_h, \bold v_h)}{||\bold v_h||_{\mathcal V^{\mathcal B}_h}} = \frac{a(\bold u, \bold v_h) - \overline a(\bold u, \bold v_h) + \overline a(\bold u, \bold v_h) - \overline{a}(\bold u_h, \bold v_h)}{||\bold v_h||_{\mathcal V^{\mathcal B}_h}} \nonumber \\
& \nonumber\\
&\mathop {\le}\limits_{\text{Holder}}& \lambda^{\bold D}_{\max}
\left(||\boldsymbol\varepsilon(\bold u) - {\overline{\boldsymbol\varepsilon(\bold u)}}||_{(L^2(\Omega))^2} +  \underbrace{||\overline{\boldsymbol\varepsilon(\bold u)} - {\overline{\boldsymbol\varepsilon(\bold u_h)}}||_{(L^2(\Omega))^2} }_{\le ||\bold u_h - \bold u||_{\mathcal V^{\mathcal B}_h}}\right ) \frac{||{\boldsymbol\varepsilon(\bold v_h)}||_{(L^2(\Omega))^2}}{||\bold v_h||_{\mathcal V^{\mathcal B}_h}}.
\end{eqnarray}
From the results (\ref{conver6}), (\ref{conver7}) and (\ref{conver8}), we get the following inequality
\begin{eqnarray}
\label{conver9}
\vspace{0.5cm}
&&\mathop{\sup}\limits_{\bold v_h \in \mathcal V^{\mathcal B}_h} \left\{
\begin{array}{l}
\vspace{0.4cm}
\lambda^{\bold D}_{\max} \left(||\boldsymbol\varepsilon(\bold u) - \overline {\boldsymbol\varepsilon(\bold u)}||_{L^2(\Omega)} + ||\bold u_h - \bold u||_{\mathcal V^{\mathcal B}_h} \right) \frac{||{\boldsymbol\varepsilon(\bold v_h)}||_{L^2(\Omega)}}{||\bold v_h||_{\mathcal V^{\mathcal B}_h}} \\
\frac{||\nabla.\bold v_h||_{L^2(\Omega)}}{{{||{{\bf{v}}_h}|{|_{\mathcal V_h^{\mathcal B}}}}}} ||p-p^{\mathcal M^*}_h||_{L^2(\Omega)} +
\frac{||\overline{\nabla. \bold v_h}||_{L^2(\Omega)}}{||\bold v_h||_{\mathcal V^{\mathcal B}_h}}||p - q_h||_{L^2(\Omega)} \end{array} \right\} \nonumber\\
&\nonumber \\
& \ge &  \beta_0 ||p_h - q_h||_{L^2(\Omega)} \ge \beta_0 (||p_h - p||_{L^2(\Omega)} - ||p - q_h||_{L^2(\Omega)})
\end{eqnarray}
for all $q_h \in \mathcal V^{**}$.\\
Thanks to the results (20), (23) in \cite{GRL1}, (79) in \cite{GRL} and the continuity property, then there exists a positive constant $\delta$ being independent on the other coefficients such that
\begin{equation}
\label{conver10}
\mathop{\sup}\limits_{\bold v_h \in \mathcal V^{\mathcal B}_h, \bold v_h \ne \bold 0} \left \{ \frac{||\nabla.(\bold v_h)||_{L^2(\Omega)}}{||\bold v_h||_{\mathcal V^{\mathcal B}_h}}, \frac{||\boldsymbol\varepsilon(\bold v_h)||_{(L^2(\Omega))^2}}{||\bold v_h||_{\mathcal V^{\mathcal B}_h}}, \frac{||\overline{\boldsymbol\varepsilon(\bold v_h)}||_{(L^2(\Omega))^2}}{||\bold v_h||_{\mathcal V^{\mathcal B}_h}}, \frac{||\bold v_h||_{(L^2(\Omega))^2}}{||\bold v_h||_{\mathcal V^{\mathcal B}_h}} \right \} \le \delta.
\end{equation}
Applying the inequality (\ref{conver10}) to two inequalities (\ref{conver3c}) and (\ref{conver9}), one obtains the following inequalities
\begin{eqnarray}
\label{conver11}
\left[ \begin{array}{l}
\vspace{0.5cm}
(\alpha_1 + \alpha_0) ||\bold w_h - \bold u||_{\mathcal V^{\mathcal B}_h}  + 2\delta ||q_h - p||_{L^2(\Omega)} + \\
\vspace{0.5cm}
 \lambda^{\bold D}_{\max} \delta ~ ||\boldsymbol{\varepsilon}(\bold u) - \overline {\boldsymbol{\varepsilon} (\bold u)}||_{(L^2(\Omega))^2} + \delta ||p^{\mathcal M^*}_h - p||_{L^2(\Omega)}\\
 \delta ||p^{\mathcal M^{**}}_h - p||_{L^2(\Omega)}
\end{array} \right] &\ge & \alpha_0 ||\bold u_h - \bold u||_{\mathcal V^{\mathcal B}_h}
\end{eqnarray}
and
\begin{eqnarray}
\label{conver12}
\vspace{0.5cm}
\frac{\alpha_0}{2} ||\overline{\boldsymbol\varepsilon(\bold u)} - {\boldsymbol\varepsilon(\bold u)}||_{L^2(\Omega)}  &+& \frac{\alpha_0}{{2\lambda^{\bold D}_{\max}}} ||p-p^{\mathcal M^*}_h||_{L^2(\Omega)} +
\frac{\alpha_0}{2 \lambda^{\bold D}_{\max}} \left(1+ \frac{\beta_0}{\delta} \right)||q_h - p||_{L^2(\Omega)} \nonumber \\
&\ge& \frac{\alpha_0\beta_0}{2\lambda^{\bold D}_{\max} \delta} ||p_h - p||_{L^2(\Omega)}  -  \frac{\alpha_0}{2}||\bold u_h -\bold u||_{\mathcal V^{\mathcal B}_h}
\end{eqnarray}
Let us the inequalities (\ref{conver11}) -(\ref{conver12}),
and use the inequality (79) in \cite{GRL}, we obtain
\begin{eqnarray}
\label{conver14}
\vspace{0.8cm}
&&(\alpha_1 + \alpha_0) \mathop {\inf }\limits_{{{\bf{w}}_h} \in \mathcal V_h^ \mathcal {B}(\lambda)} \{||\bold w_h - \bold u||_{\mathcal V_h^ \mathcal {B}}\} + \left[2\delta + {\frac{\alpha_0}{2\lambda^{\bold D}_{\max}}}{\left(1+\frac{\beta_0}{\delta} \right)}\right]\mathop {\inf }\limits_{{{p}_h} \in \mathcal V_h^{**}} \{||q_h - p||_{L^2(\Omega)} \}\nonumber \\
\vspace{0.8cm}
&+& \left(\lambda^{\bold D}_{\max} \delta + \frac{\alpha_0}{2}\right) ~ ||\boldsymbol{\varepsilon}(\bold u) - \overline {\boldsymbol{\varepsilon} (\bold u)}||_{(L^2(\Omega))^2} + \left(\delta +  \frac{\alpha_0}{{2\lambda^{\bold D}_{\max} }} \right) ||p^{\mathcal M^*}_h - p||_{L^2(\Omega)} \nonumber\\
\newline
&+& \delta ||p^{\mathcal M^{**}}_h - p||_{L^2(\Omega)} \ge  \frac{\alpha_0}{2}||\bold u_h - \bold u||_{\mathcal V^{\mathcal B}_h} + \frac{\alpha_0\beta_0}{2\lambda^{\bold D}_{\max} \delta}||p_h - p||_{L^2(\Omega)}.
\end{eqnarray}
\newline
We need to prove that there exists a positive constant $C_1$ without depending on $h$ such that
\begin{equation}
\label{conver15}
\mathop {\inf }\limits_{{{\bf{w}}_h} \in \mathcal V_h^ \mathcal {B}(\lambda)} \{||\bold w_h - \bold u||_{\mathcal V^{\mathcal B}_h}\} \le C_1 \mathop {\inf }\limits_{{{\bf{w}}_h} \in \mathcal V_h^ \mathcal {B}} \{||\bold w_h - \bold u||_{ (\mathbb H^1(\Omega))^2}\} + ||\overline{\nabla. \bold u} - \nabla. \bold u||_{L^2(\Omega)} + \frac{1}{\lambda}||p_h - p||_{L^2(\Omega)}.
\end{equation}
Let any $\bold v_h \in \mathcal V^{\mathcal B}_h$, $\bold w_h \in \mathcal V^{\mathcal B}_h(\lambda)$, we put $\bold r_h = \bold w_h - \bold v_h \in \mathcal V^{\mathcal B}_h$, then
\begin{eqnarray*}
\label{conver16}
\overline b(\bold r_h + \bold v_h,q_h) &=& b(\bold u,q_h) + \frac{1}{\lambda} \left[ \overline c(p_h,q_h) - c(p,q_h) \right ], \quad \forall q_h \in \mathcal V^{**} \nonumber \\
\Leftrightarrow \frac{\overline b(\bold r_h, q_h)}{||q_h||_{L^2(\Omega)}} &=& \frac{\overline b (\bold u- \bold v_h,q_h) + b(\bold u,q_h) - \overline b(\bold u,q_h) + \frac{1}{\lambda} \left[ \overline c(p_h,q_h) - c(p,q_h) \right ]}{||q_h||_{L^2(\Omega)}},\quad \forall q_h \in \mathcal V_h^{**}\backslash \{0\}
\end{eqnarray*}
Thanks to (\ref{phtrinh5}), we get
\begin{equation}
\label{conver17}
||\bold r_h||_{\mathcal V^{\mathcal B}_h} \le  ||\bold v_h - \bold u||_{(\mathbb H^1(\Omega))^2} + ||\overline{\nabla. \bold u} - \nabla. \bold u||_{L^2(\Omega)} + \frac{1}{\lambda}||p_h - p||_{L^2(\Omega)},
\end{equation}
which follows
\begin{eqnarray}
\label{conver18}
||\bold u - \bold w_h||_{\mathcal V^{\mathcal B}_h} &=& ||\bold u - \bold v_h - \bold r_h||_{\mathcal V^{\mathcal B}_h} \le ||\bold v_h - \bold u||_{\mathcal V^{\mathcal B}_h}  + || \bold r_h||_{\mathcal V^{\mathcal B}_h} \nonumber \\
&\le& 2 ||\bold v_h - \bold u||_{(\mathbb H^1(\Omega))^2} + ||\overline{\nabla. \bold u} - \nabla. \bold u||_{L^2(\Omega)} + \frac{1}{\lambda}||p_h - p||_{L^2(\Omega)}.
\end{eqnarray}
Hence, with $C_1 = 2$, the inequality (\ref{conver15}) is proven. 
\newline
We apply (\ref{conver15}) to (\ref{conver14}), and thank to the inequality (79) in \cite{GRL} for getting
\begin{eqnarray}
\label{conver20}
&C_2&\left[ \begin{array}{l}
\vspace{0.5cm}
2(\alpha_1 + \alpha_0) \mathop {\inf }\limits_{{{\bold{v}}_h} \in \mathcal V_h^ \mathcal {B}} \{||\bold v_h - \bold u||_{(\mathbb H^1(\Omega))^2}\} + \left[2\delta + {\frac{\alpha_0}{2\lambda^{\bold D}_{\max}}}{(1+\frac{\beta_0}{\delta})}\right]\mathop {\inf }\limits_{{{p}_h} \in \mathcal V_h^{**}} \{||q_h - p||_{L^2(\Omega)} \}\nonumber \\
\vspace{0.5cm}
+ \left(\lambda^{\bold D}_{\max} \delta + \frac{\alpha_0}{2}\right) ~ ||\boldsymbol{\varepsilon}(\bold u) - \overline {\boldsymbol{\varepsilon} (\bold u)}||_{(L^2(\Omega))^2} + \left(\delta + \frac{\alpha_0}{2\lambda^{\bold D}_{\max}} \right) ||p^{\mathcal M^*}_h - p||_{L^2(\Omega)} \nonumber\\
+ \delta ||p^{\mathcal M^{**}}_h - p||_{L^2(\Omega)} + (\alpha_1 + \alpha_0) ||\overline{\nabla. \bold u} - \nabla. \bold u||_{L^2(\Omega)}
\end{array} \right] \nonumber \\
\nonumber \\
&\ge&  ||\bold u_h - \bold u||_{\mathcal V^{\mathcal B}_h} + ||p_h - p||_{L^2(\Omega)},
\end{eqnarray}
where the positive constant $C_2$ is defined by
\begin{displaymath}
C_2 = \frac{1}{\min \left\{\frac{\alpha_0}{2}, \left( \frac{\alpha_0\beta_0}{2\lambda^{\bold D}_{\max}\delta} - \frac{\alpha_1 + \alpha_0}{\lambda}\right) \right\}},
\end{displaymath}
The coefficient $C_2$ is positive, because the Lam\'e coefficient $\lambda$ can  be chosen large enough, while $\nu$ is closed to $0.5$.\\
Besides, by using Remark $3.1$ in \cite{GRL} and the definition of $p^{\mathcal M^*}_h$ and $p^{\mathcal M^{**}}_h$, they follow
\begin{center}
$||\boldsymbol \varepsilon(\bold u) - \overline{\boldsymbol \varepsilon(\bold u)}||_{(L^2(\Omega))^2}$, $||\overline{\nabla. \bold u} - \nabla. \bold u||_{L^2(\Omega)}$, $||p^{\mathcal M^*}_h - p||_{L^2(\Omega)}$ and $||p^{\mathcal M^{**}}_h - p||_{L^2(\Omega)}$ tend to $0$,
\end{center}
as $h \to 0$.\\
Therefore, the inequality (\ref{conver1}) is proven with
\begin{displaymath}
C = C_2\max \left \{2(\alpha_0 + \alpha_1),~\left[ 2\delta + \frac{\alpha_0}{2\lambda^{\bold D}_{\max}}\left(1+\frac{\beta_0}{\delta}\right)\right]       \right\}
\end{displaymath}
and
\begin{displaymath}
\mathcal O(h) = C_2 \left[
\begin{array}{l}
\vspace{0.5cm}
\left(\lambda^{\bold D}_{\max}\delta + \frac{\alpha_0}{2}\right) ~ ||\boldsymbol{\varepsilon}(\bold u) - \overline {\boldsymbol{\varepsilon} (\bold u)}||_{(L^2(\Omega))^2} + \left(\delta + \frac{\alpha_0}{2\lambda^{\bold D}_{\max}} \right) ||p^{\mathcal M^*}_h - p||_{L^2(\Omega)} + \\
 \delta ||p^{\mathcal M^{**}}_h - p||_{L^2(\Omega)} + (\alpha_1 + \alpha_0)||\overline{\nabla. \bold u} - \nabla. \bold u||_{L^2(\Omega)}
\end{array}  \right].\quad\square
\end{displaymath}
\newline
In the following remark, we briefly recall how the scheme can be implemented for the problem (\ref{phtrinh4}) based on the displacement.\\
{\bf Remark \ref{mathproper}.2:} From (\ref{phtrinh4}$b$), we can write the pressure field $p_h = \sum\limits_{i = 1,{V_i} \in \mathcal T_h^{**}}^{{N_n}} {{p_i}{\chi _i}} $ with
\begin{eqnarray}
\label{phtrinh44a}
p_i  &=& \frac{\lambda}{m(V_i)} \int\limits_{V_i}  {(\overline { \nabla\cdot \bold{u}_h} )}  \text{d}\Omega = \frac{\lambda }{{m(V_i )}}\sum\limits_{\scriptstyle \,\,\,\,\,\,\,\,\,k = 1,\,\,\Omega _k^s  \in \mathcal T_h^*  \hfill \atop
  \scriptstyle \,\,\,\,\,\,\,\,\,\,\Omega _k^s  \cap V_i  \ne \emptyset  \hfill}^{N_s } {m(V_i  \cap \Omega _k^s )(\overline{\nabla \cdot\bold{u}_h} )\left| {_{\Omega _k^s } } \right.}.
\end{eqnarray}
Then the bilinear form $\overline b(\bold{v}_h,p_h)$ can be transformed into
\begin{equation*}
\label{phtrinh45a}
\overline b(\bold{v}_h ,p_h ) = \sum\limits_{i = 1}^{N_n } {\frac{ \lambda}{{m(V_i )}}\left( {\sum\limits_{\scriptstyle k = 1,\Omega _k^s  \in \mathcal T_h^*  \hfill \atop
  \scriptstyle \Omega _k^s  \cap V_i  \ne \emptyset  \hfill}^{N_s } {m(V_i  \cap \Omega _k^s )\left. {\left( \overline {\nabla \cdot\bold{v}_h} \right)} \right|_{\Omega _k^s } } } \right)\left( {\sum\limits_{\scriptstyle l = 1,\Omega _l^s  \in \mathcal  T_h^*  \hfill \atop
  \scriptstyle \Omega _l^s  \cap V_i  \ne \emptyset  \hfill}^{N_s } {m(V_i  \cap \Omega _l^s )\left. {\left(\overline {\nabla \cdot\bold{u}_h} \right)} \right|_{\Omega _l^s } } } \right)}.
\end{equation*}
Therefore, we arrive at a problem of finding $\bold{u}_h \in \mathcal V_h^\mathcal B$ such that
{\begin{eqnarray}
\label{phtrinh45}
\overline a(\bold{u}_h, \bold{v}_h) & +&   \sum\limits_{i = 1}^{N_n } {\frac{ \lambda}{{m(V_i )}}\left( {\sum\limits_{\scriptstyle k = 1,\Omega _k^s  \in\mathcal  T_h^*  \hfill \atop
  \scriptstyle \Omega _k^s  \cap V_i  \ne \emptyset  \hfill}^{N_s } {m(V_i  \cap \Omega _k^s )\left. {\left( \overline {\nabla .\bold u_h} \right)} \right|_{\Omega _k^s } } } \right)\left( {\sum\limits_{\scriptstyle l = 1,\Omega _l^s  \in \mathcal T_h^*  \hfill \atop
  \scriptstyle \Omega _l^s  \cap V_i  \ne \emptyset  \hfill}^{N_s } {m(V_i  \cap \Omega _l^s )\left. {\left(\overline {\nabla . \bold v_h} \right)} \right|_{\Omega _l^s } } } \right)}  \nonumber \\
 \nonumber \\
&=& (\bold{f},\bold{v}_h) \quad \quad \forall \bold v_h \in  \mathcal V^\mathcal B _h,
\end{eqnarray}}
where the solution $\bold{u}_h$ of (\ref{phtrinh45}) is the same as the solution of the problem (\ref{phtrinh4}).\\
\newline
{\bf Remark \ref{mathproper}.3:} On applying bES-FEM and bFS-FEM to linear elasticity problems, the equations can be expressed as the following linear system
\begin{equation}
\label{phtrinh46}
\left( {\begin{array}{*{20}c}
   \bold{\overline A} & {\bold{\overline B}^T }  \\
   \bold{\overline B} & { - \frac{1}{\lambda }\bold{\overline C}}  \\
\end{array}} \right)\left( {\begin{array}{*{20}c}
   {\bold{u}_h }  \\
   {p_h }  \\
\end{array}} \right) = \left( {\begin{array}{*{20}c}
   {\bold{f}_h }  \\
   0  \\
\end{array}} \right),
\end{equation}
where $\bold{\overline A}$, $\bold{\overline B}$, $\bold{\overline C}$ are matrices associated with the bilinear forms $\overline a(\cdot,\cdot)$, $\overline b(\cdot,\cdot)$ and $\overline c(\cdot,\cdot)$ respectively, and $\bold{f}_h$ is associated with the linear operator $(\bold{f},\cdot)$. This framework of bES-FEM and bFS-FEM for problems in linear elasticity has an implementation similar to that of the MINI element. However, the matrix $\bold{\overline C}$ of (\ref{phtrinh46}) is different to $\bold{C}$ in the system of linear equations 
associated with the MINI element, because the matrix $\bold{\overline C}$ of (\ref{phtrinh46}) is diagonal and each degree of freedom corresponding to the pressure can be computed by (\ref{phtrinh44a}). It follows that the matrix is positive definite.
\section{Error norms}
\label{Error}
In order to study the error and convergence of the proposed numerical methods, we introduce three error norms: the displacement error norm, the pressure error norm and the energy error norm.
\subsection{Displacement error norm}
The displacement error norm is defined by
\begin{equation}
\label{phtrinh46a}
||\bold{u} - \bold{u}_h ||_{L^2 (\Omega )}  = \left[ {\sum\limits_{T \in \mathcal T_h } {\int\limits_T {(\bold{u} - \bold{u}_h )^T (\bold{u} - \bold{u}_h )\text{d}\Omega} } } \right]^{1/2},
\end{equation}
where $\bold{u}$ is the analytical solution for the displacement and $\bold{u}_h$ is the numerical approximation.
\subsection{Pressure error norm}
The pressure error norm is written as
\begin{equation}
\label{phtrinh46b}
||p - p_h ||_{L^2 (\Omega )}  = \left[ {\sum\limits_{V \in \mathcal  T_h^{**} } {\int\limits_V {(p - p_h )^2 \text{d}\Omega} } } \right]^{1/2},
\end{equation}
where $p$ is the analytical pressure solution and  $p_h$ is the numerical solution.
\subsection{Energy error norm}
The energy error norm must take into account of the fact that some of the numerical methods solve purely for displacements but others solve additionally for pressure. The NS-FEM and ES-FEM only approximate the displacement field, hence for these two methods the evaluation of the norm follows that of \cite{GRL1} and is based on $N_s$ smoothing domains $\Omega^s_k \in \mathcal T_h^*$
\begin{equation}
\label{phtrinh46c}
||\bold{u} - \bold{u}_h ||_E  = \left\{ {\sum\limits_{\scriptstyle \hfill \atop
  \scriptstyle \Omega _k^s  \in \mathcal T_h^*  \hfill} {\int\limits_{\Omega _k^s \in T_h^*} {[\boldsymbol{\sigma}  - \boldsymbol{\overline \sigma}^{(k)}(\bold{u}_h) ]^T \,\bold{D}^{ - 1}\,[\boldsymbol{\sigma}  - \boldsymbol{\overline \sigma}^{(k)}(\bold{u}_h) ] \text{d} \Omega}} } \right\}^{1/2},
\end{equation}
where $\boldsymbol{\sigma}$ is the analytical solution for the stresses and $\boldsymbol{\overline \sigma}^{(k)}(\bold{u}_h)$, the numerical approximation to the stresses, is derived from the smoothed strain solution $\boldsymbol{\overline \varepsilon}^{(k)}(\bold{u}_h)$ defined on smoothing domains $\Omega^s_k$.\\

The MINI and bES-FEM approximate both displacement and pressure. Hence, we propose a modification to the definition of the energy error norm appropriate to each method. The norm for bES-FEM incorporates a term which depends on the pressure and is based on $N_s$ smoothing domains $\Omega^s_k \in \mathcal T_h^*$
\begin{equation}
\label{phtrinh46d}
||\bold{u} - \bold{u}_h ||_E  = \left\{ \begin{array}{l}
 2\mu \sum\limits_{\scriptstyle  \hfill \atop
  \scriptstyle \Omega _k^s  \in \mathcal T_h^*  \hfill} {\int\limits_{\Omega _k^s}{[\boldsymbol{\varepsilon}(\bold{u}) - \boldsymbol{\overline\varepsilon}^{(k)}(\bold{u}_h) ]^T \,\bold{D}\,  [\boldsymbol{\varepsilon} (\bold{u}) - \boldsymbol{\overline \varepsilon}^{(k)}(\bold{u}_h) ]~\text{d}\Omega}}  \\
\\
  + \sum\limits_{\scriptstyle  \hfill \atop
  \scriptstyle \Omega _k^s  \in \mathcal T_h^*  \hfill} {\int\limits_{\Omega _k^s } {(p(\bold{x}) - p_h ) (\nabla \cdot \bold{u}(\bold{x}) - {\left. {\left( \overline{\nabla \cdot\bold{u}_h} \right)} \right|_{\Omega _k^s } } )\,\text{d}\Omega} }  \\
 \end{array} \right\}^{1/2}.
\end{equation}
The energy error norm of the MINI method also contains a term which depends on pressure but it is evaluated on the $N_e$ triangles, $T \in \mathcal T_h$, and written as
\begin{equation}
\label{phtrinh46e}
||\bold{u} - \bold{u}_h ||_E  = \left\{ \begin{array}{l}
 2\mu \sum\limits_{\scriptstyle  \hfill \atop
  \scriptstyle T  \in \mathcal T_h  \hfill} {\int\limits_{T}[\boldsymbol{\varepsilon} (\bold{u}) - \boldsymbol{\varepsilon} (\bold{u}_h )]^T \,\bold{D}\,  [\boldsymbol{\varepsilon} (\bold{u}) - \boldsymbol{\varepsilon} (\bold{u}_h )] ~\text{d} \Omega}  \\
\\
  + \sum\limits_{\scriptstyle  \hfill \atop
  \scriptstyle T  \in \mathcal T_h  \hfill} {\int\limits_{T} {(p(\bold{x}) - p_h ) ( \nabla \cdot\bold{u}(\bold{x}) - \nabla \cdot\bold{u}_h )\,\text{d}\Omega} }  \\
 \end{array} \right\}^{1/2}.
\end{equation}
\section{Numerical results}
\label{Num}
In this section, we present some numerical results to demonstrate the efficiency and accuracy of the newly-proposed methods. For this purpose we use four benchmark problems (three cases for small deformation and a remaining one for large deformation), and compare results from bES-FEM and bFS-FEM with the results from the methods listed below.

\begin{itemize}
\item MINI - The mixed displacement-pressure finite element method with cubic bubble functions \cite{ABF}.
\item FEM - The standard FEM using three node triangular elements with linear shape functions \cite{ZT}.
\item NS-FEM - The node-based SFEM \cite{LTHK} using triangular elements.
\item ES-FEM - The edge-based SFEM \cite{LTK} using triangular elements.
{\item Q4/E4 -  The quadrilateral element implemented into four enhanced modes \cite{KT}.
\item Q4/ME2 - The mixed-enhanced formulation with
five enhanced modes. Unless otherwise noted for the results which follow the transformation matrix, $T$, used for the mixed-enhanced simulations was taken as the inverse transpose of the average Jacobian, {\it i.e.}, $T = J^{-T}_\text{avg}$ \cite{KT}.}
\item HFS-HEX8 - The hybrid finite element formulation with fundamental solutions as internal interpolation functions using linear 8-node brick elements \cite{CQY}.
\item HIS - The Hexahedral element for near-incompressibility and shear behaviour \cite{ACCF}.
\item 3D.EAS-30 - Strains are complete up to trilinear fields; the element is identical to the HR (Hellinger-Reissner) element \cite{AR}.
\item 3D.HR-18 - Hellinger-Reissner elements with the eigenvalues for 18 modes \cite{AR}.
\end{itemize}

\subsection{Cook's membrane problem}
The first benchmark test is Cook's membrane problem. This problem is often used because it serves to test how accurately a numerical method can model bending and will reveal whether or not a method is prone to volumetric locking \cite{LAMI, SR, CR,KT}. Let $\Omega$ be the convex hull
\begin{displaymath}
\Omega = \conv\{(0,0), (48,44), (48,60), (0, 44)\}.
\end{displaymath}
The domain $\Omega$ is a tapered panel (see Figure~\ref{fig:CooksMembrane}) whose left boundary is clamped, and whose right boundary is subject to an in-plane shearing load of $100$ in the $y$-direction. {Plane strain conditions are assumed. The material is described by two parameters: Young's modulus $E = 250$ and Poisson's ratio $\nu =  0.4999$. Analytical solution for this problem is not available and therefore the vertical displacement at the top conner of the right-hand boundary ({\it i.e.} the point $(48,60)$) is compared with other numerical results taken from \cite{KT}.
\begin{figure}[htbp!]
  \centering
  \caption{The domain for Cook's membrane problem, discretized with three-noded triangular elements.}
  \label{fig:CooksMembrane}
\end{figure}
\noindent A comparison between the present results and other published ones is shown Figures~\ref{fig:Cook_displacement}, \ref{fig:Cook_pressure_64} and~\ref{fig:Cook_pressure_256}.
\begin{figure}[htbp!]
  \centering
  \caption{Convergence of the displacement at point $(48,60)$ for Cook's membrane problem ($\nu =  0.4999$).}
  \label{fig:Cook_displacement}
\end{figure}
As shown in Figure~\ref{fig:Cook_displacement}, it is observed that bES-FEM can produce more accurate solution than the other methods such as MINI, ES-FEM, NS-FEM and especially mixed-enhanced strain elements \cite{KT}. ES-FEM suffers from volumetric locking. NS-FEM yields an upper bound solution and tip displacements that are oscillation-free. Unfortunately, this method cannot guarantee the stability (or the inf-sup condition addressed in Theorem 4.2) of the pressure solution which does oscillate (see Figures~\ref{fig:Cook_pressure_64},~\ref{fig:Cook_pressure_256})}.
\begin{figure}[htbp!]
  \centering
  \caption{Distribution of pressure along the line $x = 24$ for Cook's membrane problem and a mesh with 64 elements ($\nu =  0.4999$).}
  \label{fig:Cook_pressure_64}
\end{figure}
Figures~\ref{fig:Cook_pressure_64} and~\ref{fig:Cook_pressure_256} illustrate the pressure distributions through the membrane. These figures imply that the solutions of the MINI element and bES-FEM are stable, while those of ES-FEM and NS-FEM exhibit oscillations (unstable).
\begin{figure}[htbp!]
  \centering
 \caption{Distribution of pressure along the line $x = 24$ for Cook's membrane problem and a mesh with 256 elements ($\nu =  0.4999$).}
  \label{fig:Cook_pressure_256}
\end{figure}
As a further test, Cook's membrane problem is solved with distorted meshes. To generate a distorted mesh, the locations of the interior nodes of the initial mesh are modified by an irregularity factor $\text{d}$ to obtain new coordinates
\begin{eqnarray*}
x' = x + r_c\,\text{d}\,\Delta x,\\
y' = y + r_c\,\text{d}\,\Delta y,
\end{eqnarray*}
where $r_c \in [-1,1]$ is a random number; $\text{d} \in [0,0.5]$ is a distortion density; $\Delta x, \Delta y$ is the size in $x$ and $y$ directions, respectively. For two distortion densities, $\text{d} = 0.1$ and $\text{d} = 0.5$, the resulting meshes are illustrated in Figure $12$.
\begin{figure}[htbp!]
        \centering
        \begin{subfigure}[b]{6.5cm}
                \caption{$\text{d} =0.1$}
                \label{fig:distortedMeshes0.1}
        \end{subfigure}
        \begin{subfigure}[b]{6.5cm}
                \caption{$\text{d} = 0.5$}
                \label{fig:distortedMeshes0.5}
        \end{subfigure}
        \caption{Meshes of 128 four-noded triangles for Cook's membrane with two \\distortion densities. The nodes are located at the vertices and centroid of each triangle.}
        \label{fig:distortedMeshes}
\end{figure}
\begin{figure}[htbp!]
        \centering
        \begin{subfigure}[b]{6.5cm}
                \caption{Influence of mesh distortion on the accuracy of the tip displacement using 128 four-noded triangles.}
                \label{fig:influenceDistortion_disp}
        \end{subfigure}
        \begin{subfigure}[b]{6.5cm}
                \caption{Convergence of the tip displacement with a distortion density $\text{d}= 0.4$.}
                \label{fig:influenceDistortion0.4}
        \end{subfigure}
        \caption{Cook's membrane for $\nu= 0.4999$ solved on distored meshes.}
        \label{fig:influenceDistortion}
\end{figure}
\begin{figure}[htbp!]
        \centering
        \begin{subfigure}[b]{6.3cm}
                \caption{128 four-noded triangular elements}
                \label{fig:influenceDistortion_pressure_128}
        \end{subfigure}
        \begin{subfigure}[b]{6.5cm}
                \caption{4096 four-noded triangular elements.}
                \label{fig:influenceDistortion_pressure_4096}
        \end{subfigure}
        \caption{Distribution of pressure along the line $x=24$ for Cook's membrane using a mesh distortion density of $\text{d} = 0.4$.}
        \label{fig:influenceDistortion_pressure}
\end{figure}
\newline
Figure~\ref{fig:distortedMeshes} illustrates meshes consisting of 128 four-noded triangular elements generated with two distortion densities $\text{d} = 0.1$ and $\text{d} = 0.5$. The influence of irregular meshes on the displacement solution is shown in Figure~\ref{fig:influenceDistortion}. For the pressure field, it can be observed that MINI method is more sensitive to mesh distortion than bES-FEM. With the refined mesh ($64\times 64$), bES-FEM behaves well, see Figure~\ref{fig:influenceDistortion_pressure}.
\subsection{Cylindrical pipe subjected to an inner pressure}
The next benchmark problem, also considered in \cite{BB}, is a cylindrical pipe subjected to an inner pressure $p =8kN/m^2$, where its internal radius and external radius are $a = 1m$ and $b=2m$ respectively (see Figure~\ref{fig:pipe_domain}).
\begin{figure}[htbp!]
  \centering
  \caption{Model of a cylindrical pipe subjected to an inner pressure (left), and  the computational domain for this problem with symmetric conditions imposed on the left and bottom edges (right).}
  \label{fig:pipe_domain}
\end{figure}
Due to the axisymmetric nature of the problem, we only model the upper right quadrant of the pipe. We impose symmetric conditions on the left and bottom edges, the outer boundary is traction-free and a pressure is applied to the inner boundary. Plane strain conditions are applied and the Young's modulus is $E =21000 kN/m$. This problem is interesting in the nearly-incompressible case, {\it i.e.} when Poisson's ratio $\nu$ is close to $0.5$. Its domain is meshed by 3-node triangular and 4-node quadrilateral elements as shown in Figure~\ref{fig:pipe_meshes}.
\begin{figure}[htbp!]
  \centering
  \caption{Domain discretization of a cylindrical pipe subjected to an inner pressure: 256 three-noded triangular elements (left), and 128 four-noded quadrilateral elements (right).}
  \label{fig:pipe_meshes}
\end{figure}
The cylindrical pipe problem has an exact solution for the radial and tangential displacement~\cite{TiGo}
\begin{equation}
\label{phtrinh47}
\bold{u}_r (r) = \frac{{(1 + \nu)a^2 p}}{{E(b^2  - a^2 )}}\left[ {(1 - 2\nu) + \frac{{b^2 }}{r}} \right]\quad\text{and} \quad \bold{u}_{\varphi} = 0
\end{equation}
and for the stress components
\begin{equation}
\label{phtrinh48}
\boldsymbol{\sigma}_r (r) = \frac{{a^2 p}}{{b^2  - a^2 }}\left( {1 - \frac{{b^2 }}{{r^2 }}} \right),\quad\boldsymbol{\sigma} _\phi  (r) = \frac{{a^2 p}}{{b^2  - a^2 }}\left( {1 + \frac{{b^2 }}{{r^2 }}} \right),\quad \boldsymbol{\sigma}_{r\varphi }  = 0.
\end{equation}
In equations (\ref{phtrinh47})  and (\ref{phtrinh48}), $(r,\varphi)$ are the polar coordinates, and $\varphi$ is measured counter-clockwise from the
positive x-axis.\\
\newline
The rate of convergence of MINI, NS-FEM and bES-FEM is investigated for this problem and the results of this are shown in Figure~\ref{fig:pipe_norms}.
\begin{figure}[htbp!]
  \centering
       \begin{subfigure}[b]{5.42cm}
                \caption{Displacement error norms}
                \label{fig:pipe_norms_disp1}
        \end{subfigure}
        \begin{subfigure}[b]{5.42cm}
                \caption{Pressure error norms}
                \label{fig:pipe_norms_pressure1}
        \end{subfigure}
        \begin{subfigure}[b]{5.42cm}
                \caption{Energy error norms}
                \label{fig:pipe_norms_energy}
        \end{subfigure}
  \caption{Error norms of bES-FEM compared with NS-FEM and MINI method for the cylindrical pipe under the nearly-incompressible condition ($\nu =0.4999999$). The rates of convergence, \textsf{r}, can be seen in the legend of each sub-figure.}
  \label{fig:pipe_norms}
\end{figure}

According to Figures~\ref{fig:pipe_norms_disp1} and~\ref{fig:pipe_norms_pressure1}, the two convergence rates in both the displacement and the pressure error norms of bES-FEM are very high ($\ge 1.93$). The convergence rates of MINI and NS-FEM in the displacement error norm are close to 2, but their convergence rates in the pressure error norm are not as high as that of bES-FEM. Moreover, in all three norms the error in bES-FEM is lower than the error in both  MINI method and NS-FEM. Figure \ref{fig:pipe_norms_energy} confirms the convergence proof of bES-FEM as proved in Theorem 4.3.
\subsection{Nearly-incompressible block}
In this section, a nearly-incompressible block with dimensions $100 \times 100 \times 50$ is considered. The bottom face of the block is fixed and it is loaded on the top by a uniform pressure of $q = 250/$unit area, acting on an area of $20 \times 20$ at the center. By symmetry, only one quarter of the model is studied, using a tetrahedral mesh of $750$ elements with appropriate symmetry boundary conditions applied to the two interior faces. The geometry, the boundary conditions and the material parameters $E$ and $\nu$ are given in Figure~\ref{fig:block}.
\begin{figure}[htbp!]
  \centering
  \begin{subfigure}[b]{6cm}
    \centering
      \caption{Geometry and boundary conditions}
      \label{fig:block_domain}
  \end{subfigure}%
        \quad
  \begin{subfigure}[b]{6cm}
    \centering
      \caption{The mesh}
       \label{fig:block_mesh}
  \end{subfigure}
  \caption{Nearly-incompressible block.}
  \label{fig:block}
\end{figure}
The vertical displacement at the top center $P$ of the block is presented in Table~\ref{tab:block_displacements}, where the results from bFS-FEM are compared with the results from other numerical methods found in References \cite{CQY}, \cite{AR} and \cite{ACCF}.
\begin{table}[htbp!]
 \centering
    \begin{tabular}{ | l | c | }
    \hline
    numerical method & displacement \\
    \hline
     FS-FEM & 5.80E-4 \\
    \hline
    bFS-FEM & 0.02054 \\
    \hline
    HFS-HEX8 & 0.02132 \\
    \hline
    HIS      & 0.01921 \\
    \hline
    3D.EAS-30 & 0.01905\\
    \hline
    3D.HR-18 & 0.01905\\
     \hline
    \end{tabular}
\caption{Nearly-incompressible regular block, displacement at the center P of the block.}
\label{tab:block_displacements}
\end{table}
Reference \cite{ARR} reports that FS-FEM suffers from volumetric locking. In Table~\ref{tab:block_displacements} our results indicate that the bubble enrichment alleviates the locking problem. In fact, we see that bFS-FEM is softer than all but one of the other methods.
\subsection{An extension to large deformations: Case study of 2D Cook's membrane problem}
In the final test, Cook's membrane is considered for large deformations. The strain energy density of a compressible neo-Hookean material is \cite{Belytsckho1999}
\begin{equation}
\label{eq:neohookean}
\Psi\left(\mathbf{C}\right) = \dfrac{1}{2}\lambda \left(\mathrm{ln}J\right)^2-\mu\mathrm{ln}J+\dfrac{1}{2}\left(\mathrm{tr}\mathbf{C}-3\right)
\end{equation}
where $\lambda$ and $\mu$ are Lam\'{e}'s parameters as before. The bulk modulus $\kappa$ can be written in terms of these parameters: $\lambda=\kappa-\frac{2}{3}\mu$. The deformation gradient $\mathbf{F}$ is $F_{ij}=\frac{\partial x_i}{\partial X_j}$ or $\mathbf{F}=\frac{\partial\mathbf{x}}{\partial\mathbf{X}}$, and the Jacobian determinant is $J=\mathrm{det}\left(\mathbf{F}\right)$. The second Piola-Kirchhoff stress can be obtained by the first derivatives of the strain density (from equation (\ref{eq:neohookean})) with the chain rule
\begin{eqnarray}
\label{eq:pk2}
\mathbf{S}&=&2\dfrac{\partial\Psi}{\partial\mathbf{C}}=2\left(\dfrac{\partial\Psi}{\partial I_1}\dfrac{\partial I_1}{\partial\mathbf{C}}+\dfrac{\partial\Psi}{\partial I_2}\dfrac{\partial I_2}{\partial\mathbf{C}}+\dfrac{\partial\Psi}{\partial I_3}\dfrac{\partial I_3}{\partial\mathbf{C}}\right)\nonumber\\[1mm]
&=&2\left(\dfrac{\partial\Psi}{\partial I_1}+ I_1\dfrac{\partial\Psi}{\partial I_2}\right)\mathbf{I}-2 \dfrac{\partial\Psi}{\partial I_2}\mathbf{C}+2I_3\dfrac{\partial\Psi}{\partial I_3}\mathbf{C}^{-1}\nonumber\\[1mm]
&=&\mu\left(\mathbf{I}-\mathbf{C}^{-1}\right)+\lambda\mathrm{ln}J\mathbf{C}^{-1}
\end{eqnarray}
where the right Cauchy-Green deformation tensor is $\mathbf{C}=\mathbf{F}^\mathrm{T}\mathbf{F}$. The derivatives of  principal invariants with respect to the right Cauchy-Green deformation tensor $\mathbf{C}$ ($\partial I_1/\partial\mathbf{C}$, $\partial I_2/\partial\mathbf{C}$, $\partial I_3/\partial\mathbf{C}$), and the derivatives of the strain energy with respect to the principal invariants ($\partial\Psi/\partial I_1$, $\partial\Psi/\partial I_2$, $\partial\Psi/\partial I_3$) are given by \cite{Belytsckho1999}.

The elasticity tensor can be expressed in terms of the second derivatives of the strain energy density function given in equation (\ref{eq:neohookean})
\begin{equation}
\label{eq:elas_tensor}
\mathbb{C}=2\dfrac{\partial\mathbf{S}}{\partial\mathbf{C}}=4\dfrac{\partial^2\Psi}{\partial\mathbf{C}\partial\mathbf{C}}
\end{equation}
or in component form \cite{Bonnet1997}
\begin{equation}
\mathbb{C}_{ijkl}=\lambda\left(C_{ij}^{-1}C_{kl}^{-1}\right)+\left(\mu-\lambda\mathrm{ln}J\right)\left[C_{ik}^{-1}C_{jl}^{-1}+C_{il}^{-1}C_{jk}^{-1}\right]
\end{equation}
For this problem, we use the same domain $\Omega$ as in the small deformation problem with a shearing load of $1/16$ in the positive $y$-direction. The shear and bulk moduli are $\mu=0.6$ and $\kappa=1.95$, $10$, $100$, $1000$, and $10000$ respectively. Note that when the bulk modulus is $\kappa=1.95$, the neo-Hookean material is compressible, and when the bulk modulus is increased ($\kappa=10$, $100$, $1000$, and $10000$), the neo-Hookean material is approximately incompressible (Poisson's ratio is close to 0.5). The results for the proposed method are compared to the standard FEM, ES-FEM and NS-FEM with the three-noded triangular element. The numbers of elements per side are 2, 4, 8, 10, 16, 20, 32, 40, and 100 for this test.

\begin{figure}[htbp!]
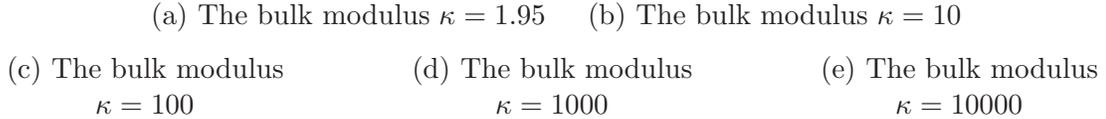

        \centering
        \begin{subfigure}[b]{5.5cm}
                \caption{The bulk modulus $\kappa=1.95$}
        \end{subfigure}\quad
        \begin{subfigure}[b]{5cm}
                \caption{The bulk modulus $\kappa=10$}
        \end{subfigure}\\

        \begin{subfigure}[b]{5cm}
                \caption{The bulk modulus $\kappa=100$}
        \end{subfigure}\quad
        \begin{subfigure}[b]{5cm}
                \caption{The bulk modulus $\kappa=1000$}
        \end{subfigure}\quad
        \begin{subfigure}[b]{5cm}
                \caption{The bulk modulus $\kappa=10000$}
        \end{subfigure}%
        \caption{Convergence of the vertical displacement at the mid-point of the right-hand boundary for Cook's membrane with the neo-Hookean model for bulk moduli ($\kappa=1.95$, $10$, $100$, $1000$, and $10000$).}
        \label{fig:cook_disp_largeDef}
\end{figure}

\begin{figure}[htbp!]
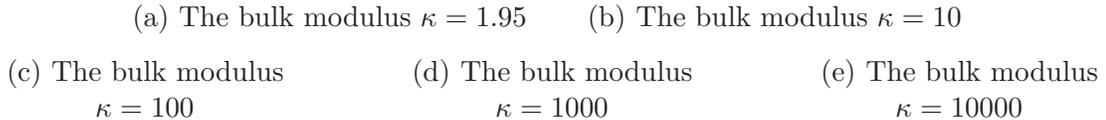

        \centering
        \begin{subfigure}[b]{5.5cm}
                \caption{The bulk modulus $\kappa=1.95$}
        \end{subfigure}\quad
        \begin{subfigure}[b]{5.5cm}
                \caption{The bulk modulus $\kappa=10$}
        \end{subfigure}\\

        \begin{subfigure}[b]{5cm}
                \caption{The bulk modulus $\kappa=100$}
        \end{subfigure}\quad
        \begin{subfigure}[b]{5cm}
                \caption{The bulk modulus $\kappa=1000$}
        \end{subfigure}\quad
        \begin{subfigure}[b]{5cm}
                \caption{The bulk modulus $\kappa=10000$}
        \end{subfigure}
        \caption{Convergence of the strain energy (logW) for Cook's membrane with the neo-Hookean model for bulk moduli ($\kappa=1.95$, $10$, $100$, $1000$, and $10000$).}
        \label{fig:cook_stne_largeDef}
\end{figure}

Figure \ref{fig:cook_disp_largeDef} illustrates the convergence of the vertical displacement at the mid-point of the right-hand boundary using both compressible and incompressible models for the proposed method, FEM, ES-FEM and NS-FEM respectively, and Figure~\ref{fig:cook_stne_largeDef} similarly illustrates the convergence of the strain energy. As shown in those figures, bES-FEM is the most robust, accurate and reliable method for both compressible and incompressible problems, compared to the conventional FEM, ES-FEM, and NS-FEM.  In the compressible problem, ES-FEM also gives relatively good convergence; however when the Poisson's ratios are close to 0.5, its convergence becomes slow. {Through the problem tested, we believe that the present method can be well applied to some relevant problems \cite{AA,BMK,NPO,PNC}.}
\section{Conclusions}\label{sec:conclusion}
We have in this paper presented the edge-based and face-based smoothed finite element methods enriched by bubble functions (bES-FEM and bFS-FEM) for nearly-incompressible elastic materials in 2D and 3D. These two methods help soften the bilinear form allowing the weakened weak ($W^2$) form to yield accurate and stable solutions. For both bES-FEM and bFS-FEM we have shown that the uniform inf-sup condition and the convergence are satisfied in the case of small deformation. {Numerical results showed, for the cases we tested, that the present method is superior to several other elements in terms of accuracy for a given number of degrees of freedom, in particular for heavily distorted meshes}.

{The proposed method is simple to implement in existing FE codes. It is efficient, and, as it does not lock even for heavily distorted triangular (simplicial) meshes which are relatively easy to generate automatically for arbitrary domains, the method is promising for incompressible problems where the structure undergoes severe deformations, as is the case during cutting and deformation of soft tissues}.

Furthermore, for problems with a curved boundary $\partial \Omega$, triangulations $\mathcal T_h$ based on simplices are not able to cover the domain $\Omega$ completely, and therefore the boundary $\partial \Omega$ is different from the boundary of $\mathcal T_h$. This issue will introduce a further error into the numerical solution. Hence, in future work, we will combine the methods presented here with NURBS functions to handle the boundary $\partial \Omega$ exactly.\\
\section*{Acknowledgements}
This  research  is  funded  by  Vietnam  National  Foundation  for  Science  and  Technology
Development  (NAFOSTED)  under  grant  number  107.02-2014.24.  The  support  is
gratefully acknowledged. The work by the senior author is partially supported by the United States NSF Grant under the award No. 1214188, and also by United States ARO contract: No.W911NF-12-1-0147. Moreover, Claire E. Heaney would like to acknowledge the financial support of EPSRC under grant EP/J01947X/1: Towards rationalised computational expense for simulating fracture over multiple scales (RationalMSFrac). These  supports  are gratefully acknowledged.


\end{document}